\documentclass[11pt,a4paper]{article}
\pdfoutput=1 
\usepackage{amsthm,amsmath,amsfonts,amssymb,graphicx,verbatim,fullpage,dsfont}
\usepackage[usenames,dvipsnames]{color}
\usepackage{enumitem}
\usepackage{ifpdf}
\ifpdf
\else
\renewcommand{\includegraphics}[1][]{\url}  
\fi

\usepackage[bookmarks,a4paper]{hyperref}
\hypersetup{%
pdftitle={Analytic properties of two-carousel systems},
pdfauthor={Ruben Bossier, Maria Vlasiou, Ivo Adan},
pdfsubject={journal paper},
pdfcreator={Maria Vlasiou},
pdfkeywords={},
pdfdisplaydoctitle=true,
bookmarks=true,%
bookmarksnumbered=true,
pdfstartview={FitH},
colorlinks=true,
citecolor=Plum,
linkcolor=BlueViolet,
urlcolor=RawSienna,
bookmarksopen=true,%
bookmarksopenlevel=1,
unicode=true,%
breaklinks=true,%
}%

\numberwithin{equation}{section}

\theoremstyle{plain}                    
\newtheorem{theorem}{Theorem}
\newtheorem{lemma}{Lemma}

\theoremstyle{remark}
\newtheorem{remark}{Remark}

\newcommand{\Dfs}[1][]{\ensuremath{F_S^{#1}}}
\newcommand{\dfs}[1][]{\ensuremath{f_S^{#1}}}  
\newcommand{\Dfa}[1][]{\ensuremath{F_A^{#1}}}
\newcommand{\dfa}[1][]{\ensuremath{f_A^{#1}}}
\newcommand{\Dfb}[1][]{\ensuremath{F_B^{#1}}}
\newcommand{\dfb}[1][]{\ensuremath{f_B^{#1}}}
\newcommand{\p}{\ensuremath{\mathbb{P}}}

\graphicspath{{graphs/}}

\begin{document}
\title{Analytic properties of two-carousel systems}
\author{Ruben Bossier$^{*}$, Maria Vlasiou$^{*,**}$, Ivo Adan$^{**,***}$}
\maketitle

\begin{tabular}{cl}
$^{*}$ &Dept.\ of Mathematics \& Computer Science, Eindhoven University of Technology\\
$^{**}$ &\textsc{Eurandom}, Eindhoven University of Technology\\
$^{***}$ &Dept.\ of Mechanical Engineering, Eindhoven University of Technology
\end{tabular}

\begin{abstract}
We present analytic results for warehouse systems involving pairs of carousels. Specifically, for various picking strategies, we show that the sojourn time of the picker satisfies an integral equation that is a contraction mapping. As a result, numerical approximations for performance measures such as the throughput of the system are extremely accurate and converge fast (e.g.\ within 5 iterations) to their real values. We present simulation results validating our results and examining more complicated strategies for pairs of carousels.
\end{abstract}

\section{Introduction}

A carousel is an automated storage and retrieval system. It consists of a number of shelves or drawers, which are linked together and are rotating in a closed loop. It is operated by a picker (human or robotic) that has a fixed position in front of the carousel. Carousel storage systems are often used to increase storage density, throughput and efficiency while reducing inventory and man-hours \cite{lebaron98}.

Carousels are widely used in flexible manufacturing systems \cite{kohler85} for storage and retrieval of small and medium-sized items \cite{seshadri94}, such as health and beauty products, repair parts of boilers for space heating, parts of vacuum cleaners and inkjet printers \cite{lebaron98}, books, shoes and many other goods. In e-commerce companies, carousels are used to store small items and manage small individual orders. An \textit{order} is defined as a set of items that must be picked together (for instance, for a single customer). Recently, carousels have been considered in healthcare applications as the means to change tools in order to enable robotic surgery without human assistance \cite{friedman07}.

Carousels are highly versatile, and come in a huge variety of configurations, sizes, and types. They can be horizontal or vertical and rotate in either only one or both directions. Although both \textit{unidirectional} (one-way rotating) and \textit{bidirectional} (two-way rotating) carousels are encountered in practice, the bidirectional types are the most common (as well as being the most efficient)~\cite{hassini03}. One of the main advantages of carousels is that, rather than having the picker travel to an item (as is the case in a warehouse where items are stored on shelves), the carousel rotates the items to the picker. While the carousel is travelling, the picker has the time to perform other tasks, such as packing or labelling the retrieved items, or serving another carousel. This practice enhances the operational efficiency of the warehouse.

Carousel models have received much attention in the literature and continue to pose interesting problems. For an extensive literature review on the topic see the surveys \cite{leduc,litvak10}. Moreover, an extensive list of references has been also made available online by Roodbergen~\cite{roodbergen_list}. One of the main topics studied in the literature is the performance of algorithms used for picking orders. Such algorithms compete with one another relative to various objectives, such as minimising the capital investment, the total motion of the carousel and the average or  worst-case time to fill an order, or maximising the utilisation of the picker and the fraction of orders that are correctly filled \cite{williams95}. For further information on picking policies, see \cite{litvak10,seshadri94,williams95} and references therein.

Studies on carousel systems offer useful insights but have a few significant drawbacks. First, most picking strategies have been studied analytically for single carousels with a dedicated picker (cf.\ the surveys \cite{leduc,litvak10} for references), while few focus on two-carousel systems that are served by a single picker \cite{hassini03,hwang91,hwang99,park03,vlasiou04,wen89}. For two-carousel systems, it is still not clear what is an optimal strategy. It is evident that the shortest picking time for the whole order on a single carousel would be achieved if the biggest gap between two consecutive items were avoided. Thus, if the carousel rotates to an item positioned at the end of the biggest gap as its starting position, the remaining distance would be minimised. However, this strategy might lead to long preparation times, compared to the time that the picker is busy at the other machine, and thus to longer waiting times for the picker, which directly translates to a lower throughput. Note that this is a potential disadvantage also for a single carousel.

Realistic systems usually involve multiple carousels (e.g.\ at the order of 10 or 20 machines) \cite{emerson81,lebaron98}. Such systems have been studied exclusively through simulation \cite{emerson81,lebaron98,williams95}. The drawbacks of this approach are evident: simulations focus on a specific configuration of a specific business, are time consuming and thus cannot be used online, and need to be repeated every time a (minor) change is considered, whether that is in the picking algorithm, the storage policy, the configuration, or the workforce allocation. Moreover, as has been observed in studies involving two carousels, insights from the analysis of picking strategies on one carousel do not immediately transfer to configurations with multiple machines \cite{mcginnis86}. Thus, systems with multiple carousels need to be studied separately. However, the complexity of the problem is such, that analytic studies beyond two carousels are intractable.

As for large groups of carousels neither analytic results nor an online (near) optimal picking algorithm exist, one hopes to use analytic insights for pairs of carousels to design a picking policy for a large-scale system. This approach makes sense since it has been empirically observed that the optimal number of carousels to assign to a single (human) picker is about two and not more than four \cite{emerson81,piasecki-IA}. Human pickers can only serve a relatively small number of machines as there is a trade-off between the utilisation of the picker (which increases as the number of carousels he serves increases since he will be rarely idle) and the throughput of the system (which decreases as the number of carousels served by a single picker increases since items will wait longer on the carousel to be picked by the picker). Thus, studying \textit{pairs} of carousels provides a realistic approximation for a system's performance measures even for large scale systems.

In this paper, we present analytic results for pairs of carousels examining a few realistic strategies. We show that carousel systems abound with contraction properties, which implies that numerical approximations for these systems are highly accurate (up to any desired level of accuracy) and converge fast (reasonable approximations can be given within 5 iterations). We study the sojourn time of the picker, namely the time needed to wait for the carousel to reach its starting position (which depends on the picking strategy considered) plus the time needed to pick a single order. For a variety of picking strategies, we show that the sojourn time distribution of the picker satisfies an integral equation that is a contraction mapping. As a result, starting from some (trivial) initial distribution, iteration of the integral equation leads to the distribution of the sojourn time, and thus to the throughput of the system. This technique is faster than simulation, more robust as it adapts easily to changes to the system's parameters, and highly accurate. Thus, it offers an alternative method to compute the performance of a pair of carousels.

The remainder of the paper is organised as follows. In the next section we introduce the problem and derive the basic integral equation we will work with. Section \ref{s:one_dir} studies a pair of unidirectional carousels in order to illustrate the steps followed and the results one can derive. Section \ref{s:bi-directional carousels} studies analytically the performance of three picking strategies for bidirectional carousels and offers some insights to their relative merits. We conclude in Section \ref{s:sim} with simulation results for an additional four strategies, which we compare for various scenarios of the order sizes.

\section{Preliminaries}
Consider two carousels operated by one picker. The picker alternates between the carousels, picking orders. The picker works at one carousel before moving to the other carousel to pick the next order (which is assumed to be always present). We consider the switching time between the carousels to be negligible. An order consists of a sequence of $n$ items which are randomly located on the carousel. Note that the number of items in an order is fixed and equal among all orders. We remove this condition in Section~\ref{s:variable}.

We are interested in the time the picker needs for one order, i.e.\ from the time he has completed picking one order at one carousel until the moment he has completed picking the next order at the other carousel. We call this the sojourn time of an order. Let $S_k$ denote the sojourn time of the $k$-th order. This time $S_{k}$ can be decomposed in two parts; namely, the time the picker might need to wait until the first item of the order is rotated by the carousel until the origin (we denote this waiting time by $W_{k}$) plus the actual picking time of the whole order (which we denote by $A_{k}$).

The waiting time of the picker is a fraction of the preparation phase of the carousel. The preparation time is the time needed from the moment an order on one carousel is picked until the moment that the carousel reaches the item that will be picked as first in the next order (on the same machine). Thus, as soon as an order is picked, the carousel starts preparing for the picking of the next order by rotating to a convenient position. The position we choose, or in other words, the item we choose to pick as first, depends on the picking strategy that is used. We denote the preparation time of the $k$-th order by $B_{k}$.

The picking time involves two actions: one is the time needed for the rotation of the carousel between the items of an order and the other is the actual time needed to remove an item from the carousel. For the order pick time $A_{k}$, we consider only the time it takes to rotate between items, ignoring the time needed to actually pick an item. In this paper, we concentrate on the travel distance (or time) and ignore the pick time; we wish to compare different strategies, and the only quantities that are modified when changing the picking strategy are the waiting time of the picker and the travel time of the carousel. The time needed to pick and package an item remains the same, and thus it can be excluded from the analysis. Thus, the time included in $A_{k}$ is the time needed to rotate to the various positions of all $n-1$ items after the first one (oftentimes referred to as ``travel time''). Alternatively, one may assume that the pick time is negligible compared to $A_{k}$ and $B_{k}$.

In this paper we examine various common picking strategies, taking into account the design of the carousel (i.e.\ whether it is able to change the direction it rotates or not). We elaborate on the particularities of each picking strategy at the relevant points later on. As in \cite{park03,vlasiou04}, we model a carousel as a circle of length 1. The carousel rotates at a unit speed; i.e., the time needed for a full rotation is also 1. Using the above notation, we can see that the sojourn time $S_{k+1}$ of the $(k+1)$-st order is the sum of the waiting time $W_{k+1}$ plus the pick time $A_{k+1}$. On the other hand, the time $W_{k+1}$ the picker needs to wait before picking the $(k+1)$-st order is given by the recursion
\begin{equation}\label{eq:W_ordered}
W_{k+1}=
(B_{k+1}-W_{k}-A_{k})^+, \qquad \mbox{for $k\geqslant 1$},
\end{equation}
where the notation $x^+$ denotes $\max\{0,x\}$ and $W_1 = B_1$; see also \cite{park03,vlasiou04}. To see this, consider that the picker needs to wait at most for a time $B_{k+1}$ for the $(k+1)$-st order until the designated carousel brings the first item to the origin.
However, while this carousel prepares (an action that started immediately upon picking the previous order on this carousel) the picker moves to the other carousel to wait $W_k$ and pick the $k$-th order in $A_k$ time. Thus, from the whole preparation time $B_k$ that the designated carousel was unavailable, a total of $W_{k}+A_{k}$ time was spent at the other carousel, and thus this time needs to be subtracted from $B_{k+1}$ as the picker was not actually waiting for the designated carousel to prepare, but was busy elsewhere. Naturally, we need the maximum of this quantity from zero as the picker cannot wait a negative time. If the designated carousel brings the first item at the origin before the picker finishes in $W_{k}+A_{k}$ time his work at the other carousel, then the first item of the $(k+1)$-st order will wait for the picker at the origin, thus incurring a zero waiting time for that order.

%
It is important to note here that \eqref{eq:W_ordered} has not been studied before in this form, despite its evident relation to the equations studied in \cite{park03,vlasiou} and related work. The main difference is though that the right-hand side of Equation \eqref{eq:W_ordered} involves three random variables that in previous studies were assumed to be independent. In this work, $A_k$ depends on $W_k$ that appears in the right-hand side of Equation \eqref{eq:W_ordered}. The form of the dependence is unknown, leading to added complications compared to previous work on this Lindley-type equation. In order to overcome the difficulties posed by this unknown form of the dependence, we translate \eqref{eq:W_ordered} into an equation for the sojourn time $S_{k+1}$ of the picker. Keeping in mind that $S_k=W_k+A_k$, we now see that for $k\geqslant1$
\begin{equation}\label{main_ordered}
S_{k+1} = (B_{k+1}-S_k)^++A_{k+1},
\end{equation}
where $S_1 = B_1+A_1$. Although the random variables $A_{k+1}$ and $B_{k+1}$ that appear on the right-hand side of \eqref{main_ordered} depend on one another, their dependence is known (and is determined by the picking strategy followed), and moreover, they are both independent of the sojourn time $S_k$ appearing on the right-hand side of \eqref{main_ordered}, thus making it possible to analyse this equation.
In fact, we wish to derive the limiting distribution of $S_k$ as $k$ tends to infinity, which can be shown (by standard regenerative arguments) to exist and to be unique; for a proof of this statement in all generality see \cite{vlasiou07a}. In steady state, \eqref{main_ordered} becomes
\begin{equation}\label{main}
S\stackrel{\mathcal{D}}{=}(B-S)^++A,
\end{equation}
where the generic random variable $S$ is distributed according to the limiting distribution of the sojourn time. This is the main equation we study here. Although the focus is on the steady-state equation \eqref{main}, it should be noted that the approach developed in this paper to solve \eqref{main} also gives us the transient distribution of $S_k$, see Remark \ref{rem:transient}.

The sojourn time is an important performance measure as its reciprocal gives the throughput of the system. The preparation time $B$ and the pick time $A$ depend on the picking strategy followed. In the sequel, we study various picking strategies. For each strategy, we derive the distributions for $A$ and $B$, and by using \eqref{main}, we derive an integral equation that the sojourn time satisfies (which naturally is different for each case). We then study the properties of that integral equation and see how we can use this information to derive practical numerical approximations of the throughput that are both fast and accurate.

We use the following notation. By \Dfs\ we denote the distribution function of the random variable $S$, while its density is given by \dfs. We denote the distributions and densities of other random variables in a similar fashion. We then have that $\Dfs(0)=0$ and for $0\leqslant x \leqslant 1$
\begin{align} \label{eq2}
\Dfs(x)&=\p[(B-S)^++A\leqslant x]=\int_{0}^{x}\p[(B-S)^+\leqslant x-y|A=y]\dfa(y)dy\nonumber\\
&=\int_{0}^{x}\p[B-S\leqslant x-y|A=y]\dfa(y)dy\nonumber\\
&=\int_{0}^{1}\int_{0}^{x}\p[B\leqslant x-y+z|A=y]\dfa(y)dy \dfs(z)dz.
\end{align}

We assume that a picking policy does not lead to a distance more than a full circle. In other words, we assume that $0 \leqslant S \leqslant  1$, and thus $\Dfs(1)=1$. Should a picking policy lead to such a situation, we simply choose to make one full rotation.

\section{Unidirectional carousels}\label{s:one_dir}

In this section, we study the simplest case possible and under the simplest strategy possible. The reasons for doing so are as follows: as will be evident in the sequel, the computations are cumbersome and thus we focus on the easiest case to help the reader along; additionally, the insights gained from this case carry on also to more complicated strategies and bidirectional carousels. We make these connections evident later on. Thus presenting the simple case builds up the intuition for later cases. 

For unidirectional carousels, one can contrive several strategies that position the carousel at a convenient starting point. E.g., one could choose to position the carousel at the optimal point for a single order and a single machine (i.e.\ the end of the biggest gap); however, this strategy is not necessarily optimal for a \textit{pair} of carousels as it may lead to a decrease in throughput. The time needed to rotate the carousel to the end of the biggest gap might be comparatively too long. In this section, we study the simplest strategy possible. Namely, we consider the case where the preparation phase is simply to have the carousel rotate to the first item that it encounters and wait there for the picking to commence.

\subsection{Integral equation}
For orders of size 1, the problem under this setting (a pair of unidirectional carousels rotating to the first item on their way for the preparation phase) has been studied in \cite{park03}. Thus, in the remainder of this section we consider the case of $n\geqslant 2$. One can readily see that for $n\geqslant 2$ we have that
\begin{equation}\label{eq3}
\dfa(y)=\begin{cases}
n(n-1)(1-y)y^{n-2}, & \mbox{for $0<y<1$};\\
0, & \mbox{otherwise}.
\end{cases}
\end{equation}
We see that this is the density of the sample range of $n$ ordered numbers; there are $n$ possible items to start from and thus $n-1$ possible items to finish picking the order. The probability that the other $n-2$ items are located on one side between these two tagged items is $y^{n-2}$. The probability that the starting point is located on the other side is $1-y$.

As is also evident, we have that
\begin{equation}\label{eq:1dir_CondB|A}
\p[B\leqslant x|A=y]=\frac{x}{1-y}, \quad 0\leqslant x\leqslant 1-y.
\end{equation}

Substituting these two expressions into \eqref{eq2} we have that for $0\leqslant x\leqslant 1$,
\begin{align}\label{eq1}
\nonumber \Dfs(x)&=\int_0^{1-x}\int_0^x(x-y+z)n(n-1)y^{n-2}dy \dfs(z)dz+\int_{1-x}^1\int_0^xn(n-1)(1-y)y^{n-2}dy \dfs(z)dz\\
\nonumber &=x^n+nx^{n-1}\int_0^{1-x}z \dfs(z)dz+nx^{n-1}(1-x)\int_{1-x}^1 \dfs(z)dz\\
\nonumber &=x^n+nx^{n-1}\left[(1-x)\Dfs(1-x)-\int_0^{1-x}\Dfs(z)dz\right]+nx^{n-1}(1-x)[1-\Dfs(1-x)]\\
&=nx^{n-1}-(n-1)x^n-nx^{n-1}\int_0^{1-x}\Dfs(z)dz.
\end{align}

\subsection{Contraction mapping}\label{ss:1dir_contr}
Consider the Lebesgue space $\mathcal{L}=\mathcal{L}^\infty([0,1])$ of all measurable and bounded functions on $[0,1]$ endowed with the  sup norm. We can prove the following.

\begin{theorem}\label{th:1}
For $n \geqslant 2$ and $0 \leqslant x \leqslant 1$, the operator $\Omega:\mathcal{L}\rightarrow\mathcal{L}$ that is defined by
\begin{equation*}
\Omega F(x)=nx^{n-1}-(n-1)x^n-nx^{n-1}\int_0^{1-x}F(z)dz,\quad \mbox{with $F\in \mathcal{L}$},
\end{equation*}
is a contraction mapping on $\mathcal{L}$ with contraction constant given by $(\frac{n-1}{n})^{n-1}$.
\end{theorem}

\begin{proof}
Let $F_1(x)$ and $F_2(x)$ be two functions in $\mathcal{L}$ with $|F_1(x)-F_2(x)|\leqslant \delta$ for all $x\in [0,1]$. We show that for some constant $c<1$,
$$
||\Omega F_1-\Omega F_2||\leqslant c\ ||F_1-F_2||,
$$
where
$$
||F|| = \sup_{x\in [0,1]} |F(x)|.
$$
In order to compute $||\Omega F_1-\Omega F_2||=\sup_x |\Omega F_1(x)-\Omega F_2(x)|$, observe that
\begin{align*}
|\Omega F_1(x)-\Omega F_2(x)|&\leqslant |nx^{n-1}\int_0^{1-x}F_1(z)-F_2(z)dz|\leqslant nx^{n-1}\int_0^{1-x}|F_1(z)-F_2(z)|dz\\
&\leqslant nx^{n-1}\int_0^{1-x}\delta dz=\delta nx^{n-1}(1-x).
\end{align*}
The function $g(x)=nx^{n-1}(1-x)$ for $n\geqslant 2$ is positive for all $x\in (0,1)$ and continuous with its maximum value achieved at $x=\frac{n-1}{n}$ and given by $g(\tfrac{n-1}{n})=(\frac{n-1}{n})^{n-1} < 1$. Thus $\Omega$ is a contraction mapping on $\mathcal{L}$ with the constant $(\frac{n-1}{n})^{n-1}$.
\end{proof}

\begin{remark}\label{rem:extendSpace}
In case the sojourn time of the picker is not bounded by 1, one needs to extend all proofs to the Lebesgue space $\mathcal{L}^\infty(\mathds{R})$, and define appropriately the opearator $\Omega$. This is naturally possible also in the case considered in Theorem~\ref{th:1}, but unnecessary. In this case, since \Dfs\ is the distribution of a random variable bounded by 1, one would define $\Omega$ to be equal to 1 for all $x>1$ and equal to 0 for all $x<0$. Then, the absolute difference of $|F_1(x)-F_2(x)|$, and thus also the norm $||\Omega F_1-\Omega F_2||$, is equal to zero outside $[0,1]$.
\end{remark}

\begin{figure}
\begin{minipage}{\textwidth}
\includegraphics[width=0.49\textwidth]{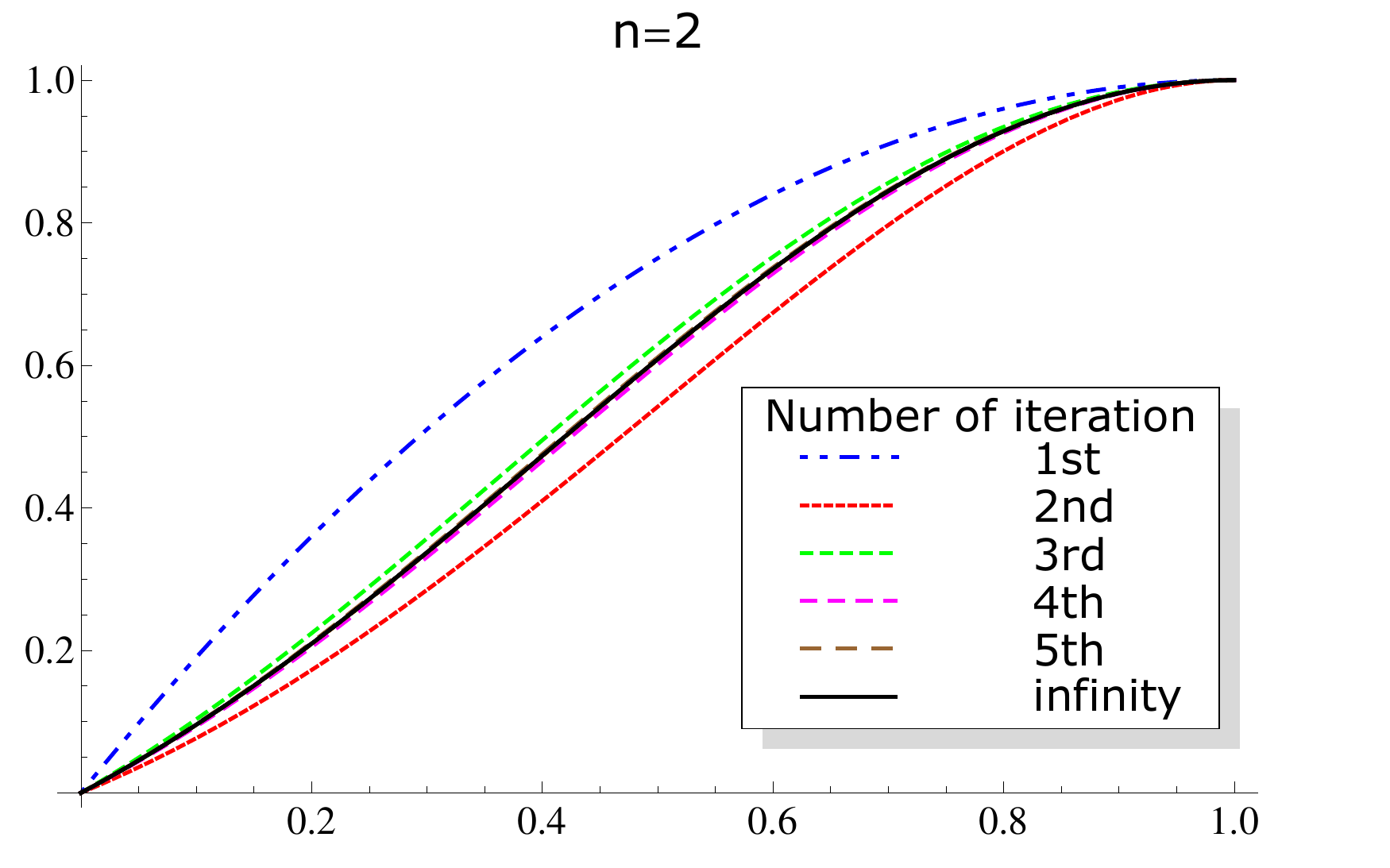}\hfill
\includegraphics[width=0.49\textwidth]{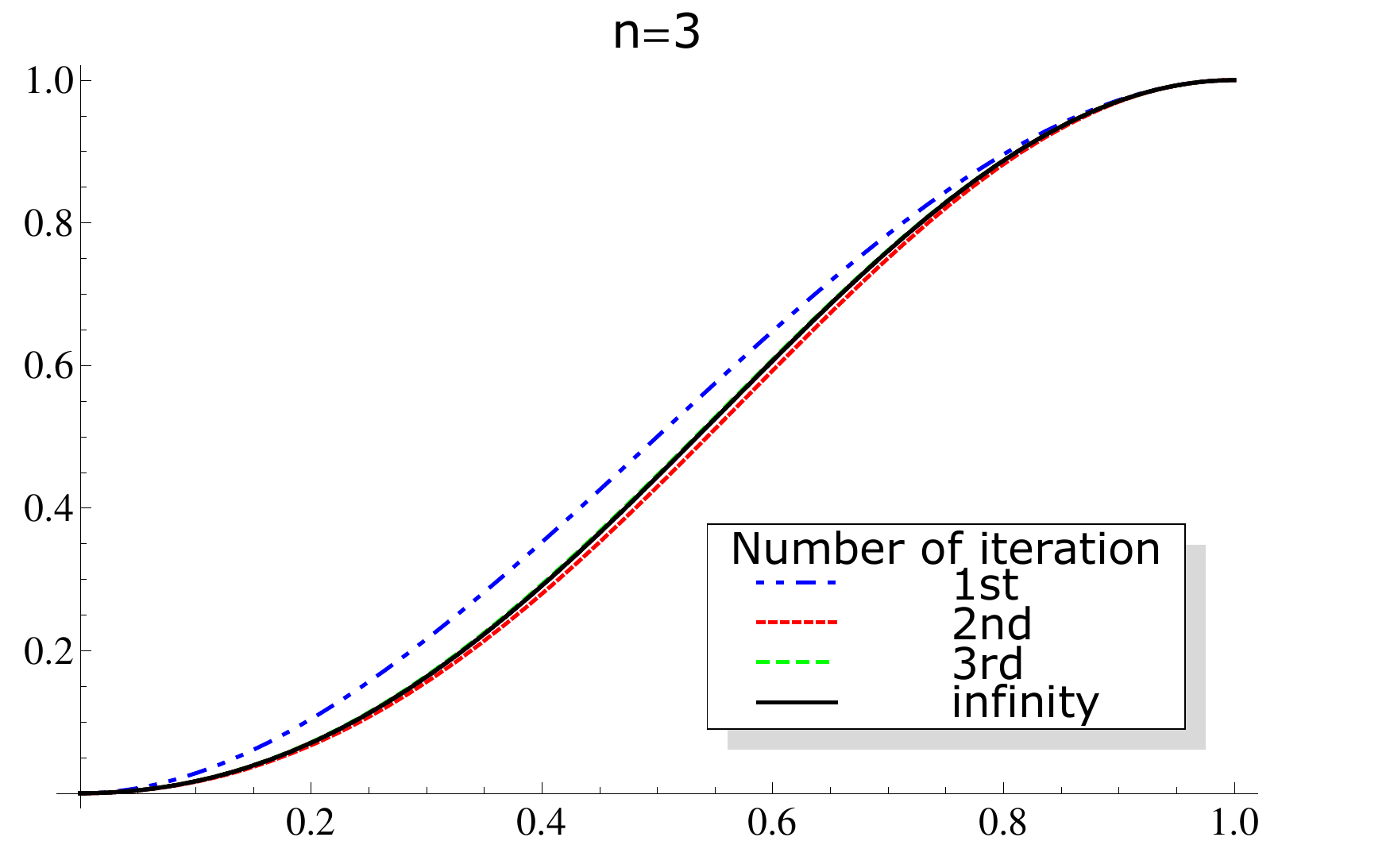}
\end{minipage}
\begin{minipage}{\textwidth}
\includegraphics[width=0.49\textwidth]{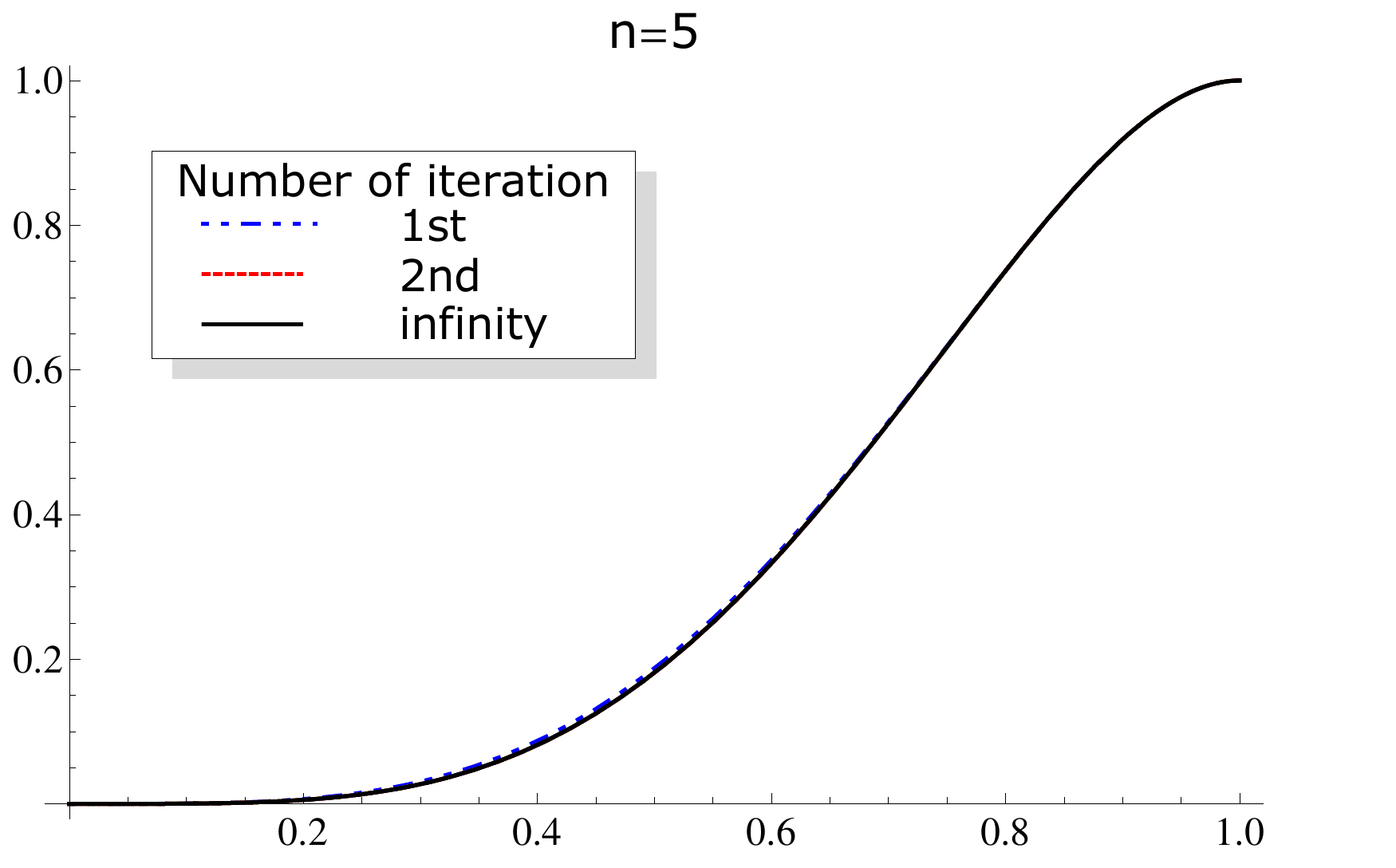}\hfill
\includegraphics[width=0.49\textwidth]{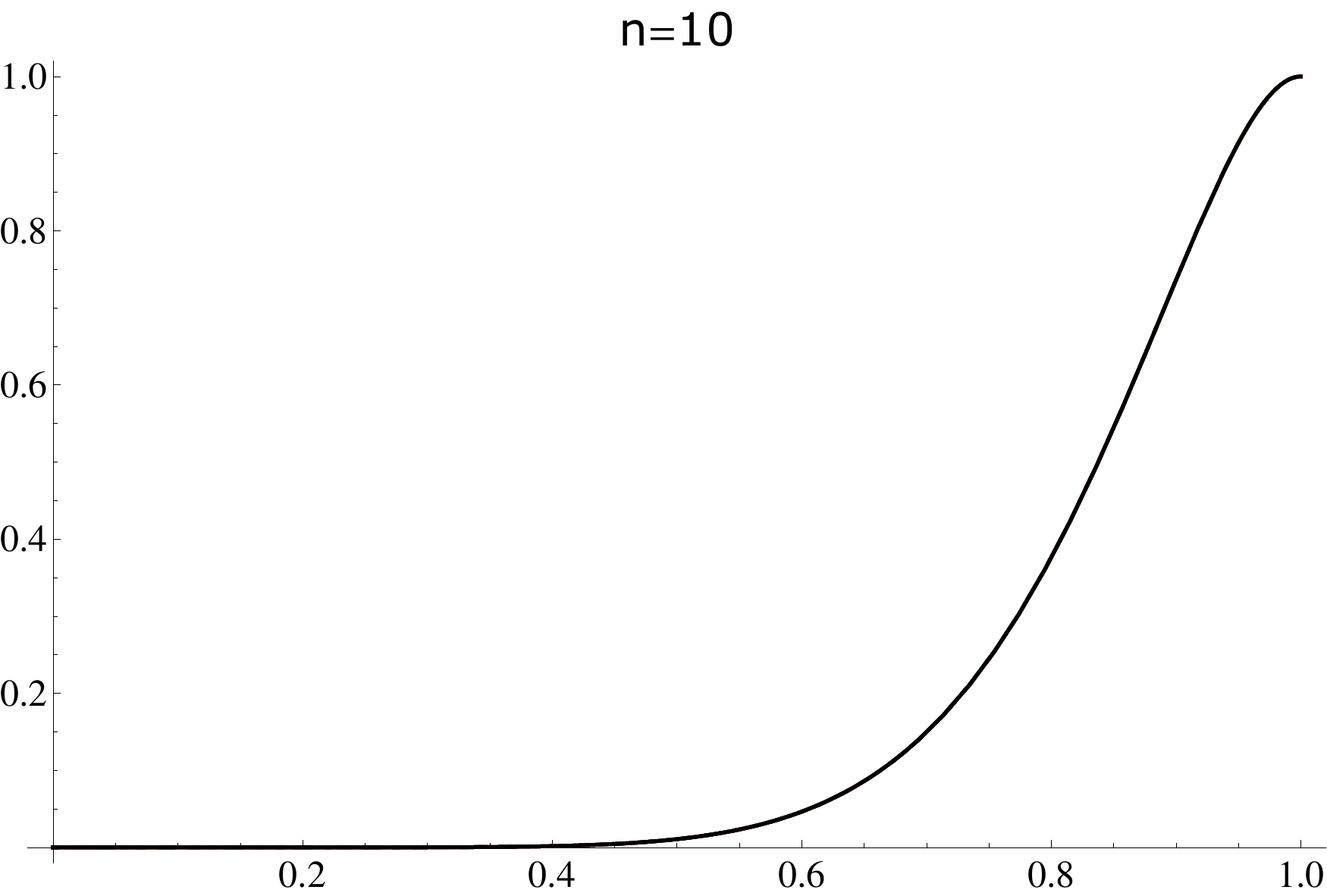}
\end{minipage}
\caption{From left to right and top to bottom we see the convergence for the values $n=2$, $n=3$, $n=5$, and $n=10$. The iterations were initiated with the upper limit $F_1(x)=nx^{n-1}-(n-1)x^n$.\label{fig2}}
\end{figure}

Theorem~\ref{th:1} is useful in determining the sojourn time distribution \Dfs\ by successive substitutions. Namely, by starting with some initial function in $\mathcal{L}$, e.g.\ with $F_0(x)=0$, one can obtain successive approximations of \Dfs; see also Remark~\ref{rem:specialF}. Observe from \eqref{eq1} that $0\leqslant \Dfs(x)\leqslant nx^{n-1}-(n-1)x^n$ and that $\lim_{n \to \infty} (\frac{n-1}{n})^{n-1}=\frac{1}{e} \simeq 0.367879$. Thus iteration converges fast. Starting with $F_0(x)=0$, we find $F_1(x)=nx^{n-1}-(n-1)x^n$ (the upper limit); iteration produces alternatingly a function that is bigger or smaller than the limit. In Figure~\ref{fig2} we see the convergence to the limiting function for various values of the number of items $n$.

The Lebesgue space $\mathcal{L}$ is a Banach space. Consequently, by the Banach Fixed Point Theorem \cite{istratescu-FPT}, we have that there is a unique measurable bounded function $F$ in $\mathcal{L}$ that is the fixed point of $\Omega$. The set of continuous and bounded functions on $[0,1]$ with the sup norm is also a Banach space, since it is a closed subspace of $\mathcal{L}$. Since \Dfs\ is a fixed point of $\Omega$, as it can be shown that there exists a unique limiting sojourn time distribution, and since it is c\`{a}dl\`{a}g and thus in $\mathcal{L}$, we conclude that the fixed point of  $\Omega$ that we converge to by iteration is necessarily the steady-state sojourn distribution. This argument comprises the main idea of all subsequent results. For each case, we examine if the related mapping is a contraction, and if so, since $\Dfs$ is in $\mathcal{L}$, then iteration leads to it.

\begin{remark}\label{rem:transient}
Observe that this iterative approach gives us the transient distribution of the sojourn time $S_k$ of the $k$-th order for a given distribution for $S_1$. Additionally note that we have also computed a geometric upper bound for the rate of convergence to the invariant distribution, namely the constant of the contraction mapping $(\frac{n-1}{n})^{n-1}$, which justifies our claim that this numerical procedure is numerically accurate and converges (geometrically) fast to the steady-state distribution.
\end{remark}
\begin{remark}\label{rem:interpretTerm}
Equation \eqref{eq1} can be interpreted as follows. Consider \eqref{eq3} and observe that \eqref{eq1} can be written as
\begin{equation*}
\Dfs(x)=\Dfa(x)-nx^{n-1}\int_0^{1-x}\Dfs(z)dz.
\end{equation*}
Moreover, we see from \eqref{main} that
\begin{equation*}
\p[S\leqslant x]=\p[(B-S)^++A\leqslant x]=\p[A\leqslant x]-\p[(A\leqslant x)\cap((B-S)^++A> x)],
\end{equation*}
and thus
\begin{equation*}
\p[(A\leqslant x)\cap((B-S)^++A> x)]=nx^{n-1}\int_0^{1-x}F(z)dz.
\end{equation*}
\end{remark}

\subsection{Variable number of items}\label{s:variable}

So far we have assumed that an order consists of a fixed number of items equal to $n$. Here we show that the analysis follows also in the case that the size of an order is random, equal to $N$, and has some known distribution. Specifically, assume that the probability $p_m=\p[N=m]$ that an order has $m$ items is known. Denote by \Dfs\ the (arbitrary) sojourn time distribution of the picker, and by $\Dfs[(m)]$ the sojourn time distribution conditioned on the event that an order has size $m$. Thus, $\Dfs(x)=\sum_m p_m \Dfs[(m)](x)$.

As before, \eqref{eq:1dir_CondB|A} holds, and analogously to \eqref{eq2}, we have that
\begin{equation*}
\Dfs[(m)](x)=\int_0^1\int_0^x\p[B\leqslant x-y+z|A=y)\dfa(y)\dfs(z)dydz.
\end{equation*}
In the special case that an order consists of precisely one item, we have that the preparation time has the uniform distribution on $[0,1]$ and that the travel time is equal to zero. Thus,
\begin{equation*}
\Dfs[(1)](x)=\Dfb(x) F(0)+\int_0^1\p[B\leqslant x+z]\dfs(z)dz.
\end{equation*}
Expressing \dfs\ as a sum of conditional densities, we can derive an integral equation for $\Dfs[(m)]$ along the same lines as for \eqref{eq1}. We then get that for $m=1$,
\begin{equation}\label{eq:1dir_var_1}
\Dfs[(1)](x)=1-\sum_i p_i\int_0^{1-x} \Dfs[(i)](z) dz,
\end{equation}
and for $m\geqslant 2$,
\begin{equation}\label{eq:1dir_var_2}
\Dfs[(m)](x)= m x^{m-1}-(m-1)x^m-mx^{m-1}\sum_i p_i \int_0^{1-x} \Dfs[(i)](z) dz.
\end{equation}

From Theorem~\ref{th:1} it immediately follows that for all fixed $m$,  $\Dfs[(m)]$ satisfies a contraction mapping and can thus be approximated by successive iterations. Thus, in case the number of items is not fixed, but has a discrete distribution, the sojourn time distribution \Dfs\ can be estimated by first estimating the conditional distributions $\Dfs[(m)]$ from the system of equations formed by \eqref{eq:1dir_var_1} and \eqref{eq:1dir_var_2} by iterating the vector $(\Dfs[(1)], \Dfs[(2)], \ldots, \Dfs[(m)])$ and then multiplying appropriately with the distribution of the size order to derive the unconditional sojourn time. Similar computations can be made also for the cases examined in the following section.

\begin{remark}
From \eqref{eq:1dir_var_1} and \eqref{eq:1dir_var_2} we can readily observe the following interesting relations:
\begin{align*}
\Dfs[(1)](0)&=1-\int_0^{1} \Dfs(z)dz=\int_0^1 z \dfs(z)dz = E[S],\\
\Dfs[(1)'](0)&=\sum_i p_i \Dfs[(i)](1)=1,\\
\Dfs[(1)'](1)&=\sum_i p_i \Dfs[(i)](0)= p_1 \Dfs[(1)](0)= p_1 E[S],\\
\Dfs[(2)'](0)&=2(1-\sum_i p_i \int_0^1 \Dfs[(i)](z) dz)=2\left(1-\int_0^1 \Dfs(z)dz\right) = 2 E[S],\\
\Dfs[(2)'](1)&=2(1-1-0)+2\sum_i p_i \Dfs[(i)](0) = 2 p_1 E[S].
\end{align*}
\end{remark}

\section{Bidirectional carousels}\label{s:bi-directional carousels}
Most carousels nowadays are able to rotate in two directions \cite{piasecki-IA}. This added freedom leads to a variety of strategies that can be followed in order to pick an order efficiently. In the sequel, we examine three reasonable strategies that can be used for bidirectional carousels.
\begin{enumerate}[topsep=0pt, itemsep=0pt]
\item \textsl{Go to the nearest item:} In the first strategy, the carousel moves to the nearest item from the two that lie on the left and on the right of the starting point at the beginning of the preparation time, and after picking the first item, it proceeds by rotating in the same direction in order to collect the rest of the order. Thus no direction reversals occur within the picking of one order.
\item \textsl{Shortest single direction after the nearest item:} In the second strategy, the carousel also moves to the nearest item first. However, from that point, we decide whether collecting all items clockwise or collecting them all counter-clockwise is shorter and we choose the shortest route between these two. At that point a reversal of the original rotation direction might occur (recall that for a single carousel it is known that the optimal strategy admits at most one turn \cite{bartholdi86}; although this fact is not proven for pairs of carousels, it seems rather evident).
\item \textsl{Avoid the largest gap:} The last strategy we present is to position the carousel at a point of the largest gap. Evidently, this strategy leads to the shortest travel time per order, but is not necessarily optimal (as positioning the next carousel to the ``optimal'' starting point might lead to long preparation times, thus wasted time for the picker, and thus lower throughput).
\end{enumerate}
All three strategies are equivalent if orders consist of precisely one item.  With minor modifications, the analysis in this case follows \cite{park03}. Thus, in the sequel we assume that $n \geqslant 2$. To facilitate the flow and emphasise the main points, we present most proofs in the appendix.

\subsection{Single direction: go to the nearest item}\label{s:NI}
In the first strategy we examine the scenario in which the carousel rotates in one direction: i.e.\ the direction towards the item nearest to the origin (at the beginning of the preparation time). As we rotate in a single direction (for each order), we have that \eqref{eq3} for \dfa\ holds and that
\begin{equation}\label{eq:condB}
\p[B\leqslant x|A=y]=\begin{cases}
\frac{2x}{1-y}, &0\leqslant x\leqslant \frac{1-y}{2}\\
1, &x> \frac{1-y}{2};
\end{cases}
\end{equation}
i.e., that the preparation time, conditioned on the event that the travel time is equal to $y$ is uniformly distributed on $[0,\frac{1-y}{2}]$. Note that the domain is half of the remaining distance, as we move to the nearest item first (which cannot be at a distance more than $\frac{1-y}{2}$ units from the origin).

By substituting these two equations into \eqref{eq2} we can now readily prove the following statement.
\begin{lemma}\label{lem:1}
For $n \geqslant 2$ and $0 \leqslant x \leqslant 1$, the sojourn time of the picker satisfies the following integral equation:
\begin{equation}\label{eq:int}
\Dfs(x)=nx^{n-1}-(n-1)x^n-2n(n-1)\int_{(2x-1)^+}^{x}y^{n-2}\int_{0}^{\frac{1}{2}+\frac{y}{2}-x}\Dfs(z)dzdy,
\end{equation}
\end{lemma}
\begin{proof}
  See Appendix~\ref{app:1}.
\end{proof}
As before, the metric space $\mathcal{L}$ with the sup norm is complete. Starting with \eqref{eq:int}, we can prove the following.

\begin{theorem}\label{th:2}
The operator $\Omega:\mathcal{L}\rightarrow\mathcal{L}$ that is defined by
\begin{equation*}
\Omega F(x)=nx^{n-1}-(n-1)x^n-2n(n-1)\int_{(2x-1)^+}^{x}y^{n-2}\int_0^{\frac{1}{2}+\frac{y}{2}-x}F(z)dzdy,\quad
\mbox{with $F\in\ \mathcal{L}$}
\end{equation*}
for $0 \leqslant x \leqslant 1$ is a contraction mapping on $\mathcal{L}$ with contraction constant $\left(\frac{n-1}{n}\right)^{n-1}$.
\end{theorem}
\begin{proof}
  See Appendix~\ref{app:2}.
\end{proof}

\begin{figure}
\begin{minipage}{\textwidth}
\includegraphics[width=0.49\textwidth]{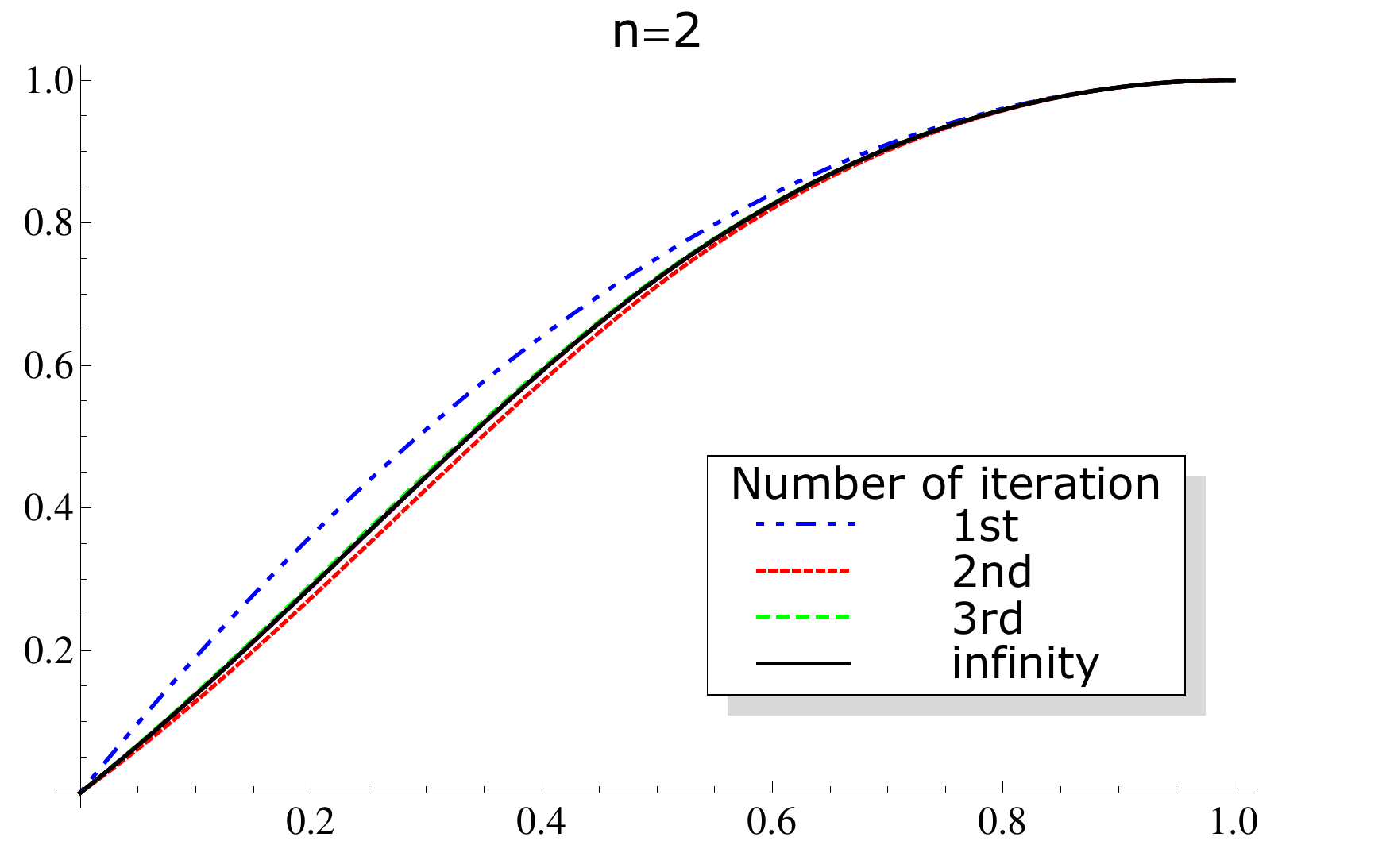}\hfill
\includegraphics[width=0.49\textwidth]{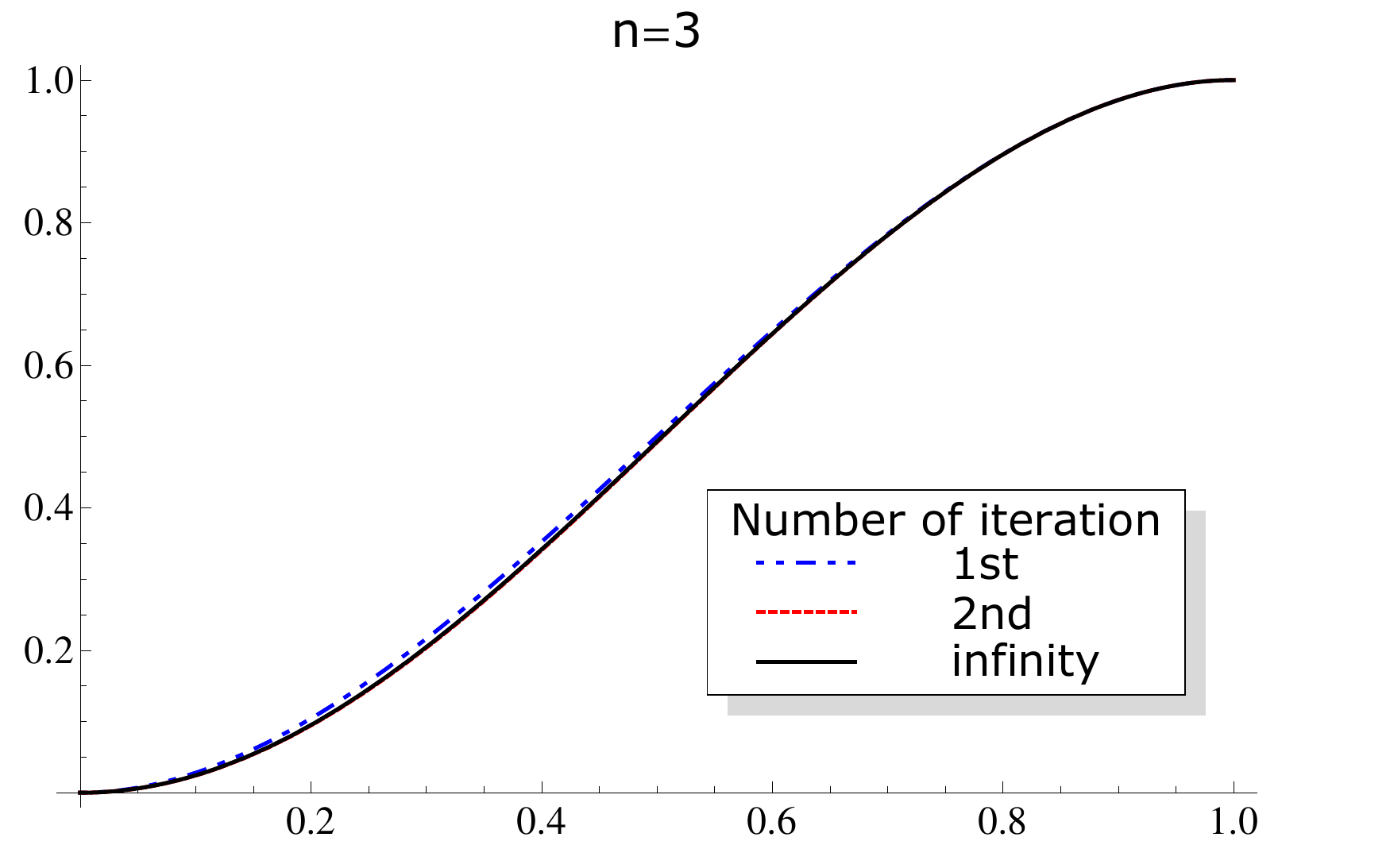}
\end{minipage}
\begin{minipage}{\textwidth}
\includegraphics[width=0.49\textwidth]{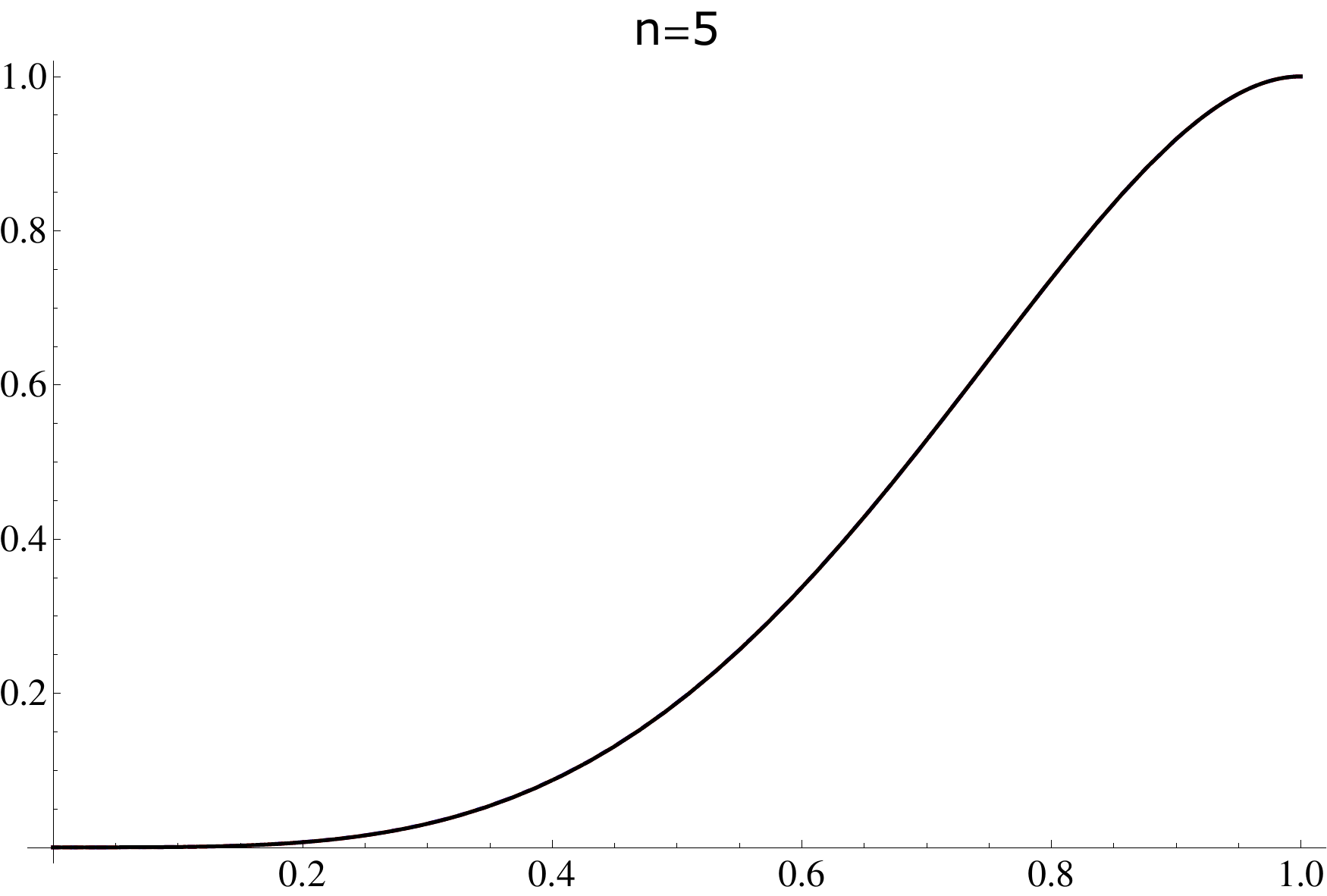}\hfill
\includegraphics[width=0.49\textwidth]{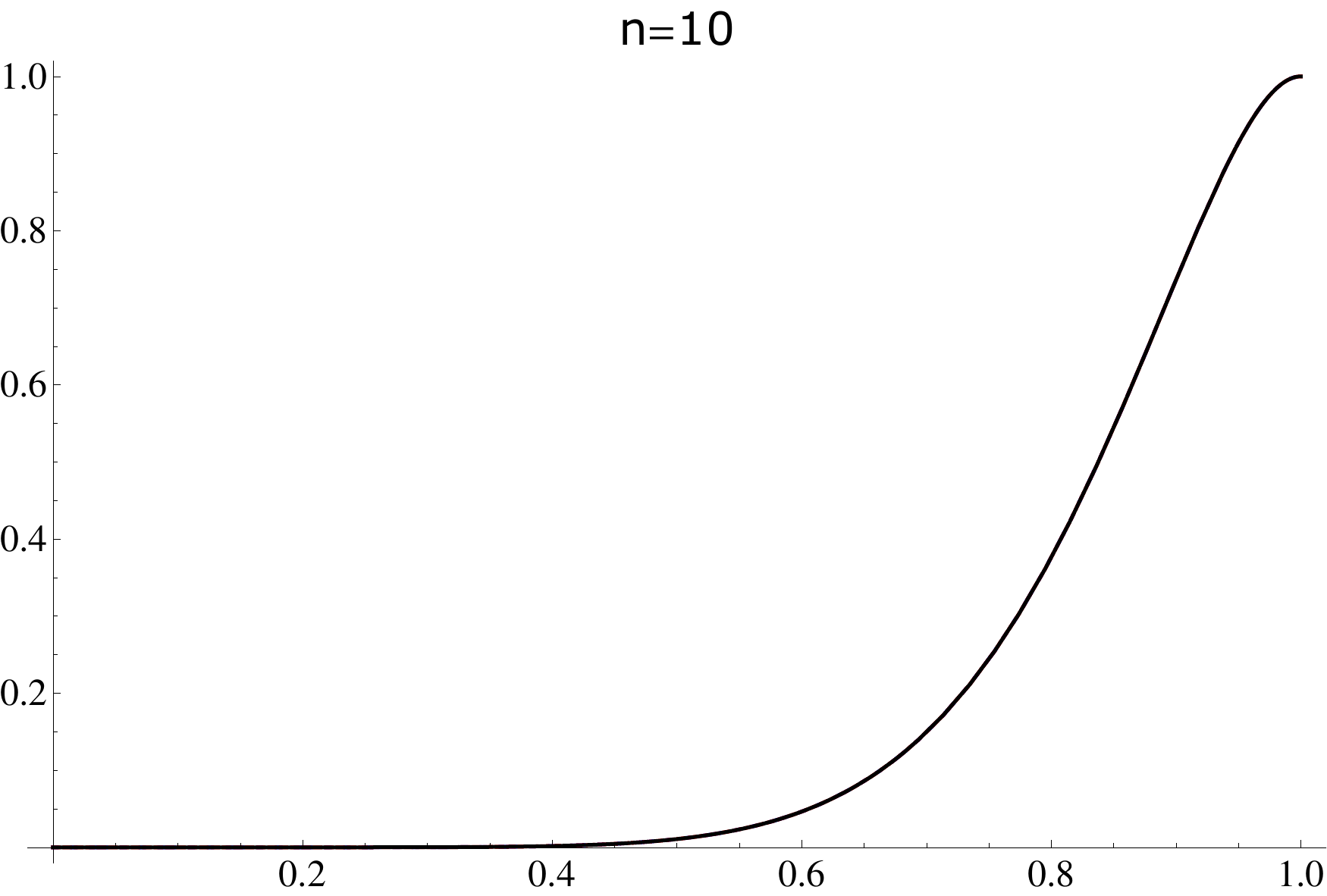}
\end{minipage}
\caption{From left to right and top to bottom we see the convergence for the values $n=2$, $n=3$, $n=5$, and $n=10$. The iterations were initiated with the upper limit $F_1(x)=nx^{n-1}-(n-1)x^n$.\label{fig4}}
\end{figure}

From \eqref{eq:int} we see that $0\leqslant \Dfs(x)\leqslant nx^{n-1}-(n-1)x^n$. Iteration will result in an alternating upper bound and lower bound for  \Dfs. For example, assume that the picking of the first order starts immediately with collecting items, i.e.\ that the preparation is already done. For the strategy considered in this section, this implies that the carousel has already moved to the nearest item before we start time. Then, we have that for $0\leqslant x <1$, $F_0(x)=0$, and we find the first estimate $F_1(x)=nx^{n-1}-(n-1)x^n$. In Figure~\ref{fig4} we see the evolution of the iteration for various values of $n$.
\begin{remark}\label{rem:specialF}
Note that we can start with any function $F_0$ that is in $\mathcal{L}$, even if our choice has no physical interpretation. Due to the contraction principle, we will converge to the steady-state distribution \Dfs, although in such a case we lose the interpretation that the $i$-th iteration yields the transient distribution for the $i$-th sojourn time. Thus, a special function $F_0$ is not required, but is convenient for an interpretation in terms of transient analysis.
\end{remark}

\begin{remark}
From the figure, we can make the following interesting observation. The convergence for this strategy is faster than what was observed for unidirectional carousels. The difference between these two cases lies only on the time it takes to rotate to the first item (i.e.\ the preparation time $B$). Recall that $\p[S\leqslant x]=\p[(B-S)^++A\leqslant x]$ and that for both cases we started with the first estimate $F_1(x)=nx^{n-1}-(n-1)x^n$ which is equal to $\Dfa(x)$. Since though for this case we have that the preparation time $B$ is smaller than that of the first case (for unidirectional carousels), we have that $(B-S)^+$ is closer to zero and thus \Dfa\ is a better estimate in the bidirectional case -- which explains why the convergence is faster.
\end{remark}

\begin{remark}
As in Remark~\ref{rem:interpretTerm}, we can write
\begin{equation*}
\Dfs(x)=\Dfa(x)-2n(n-1)\int_{(2x-1)^+}^{x}y^{n-2}\int_0^{\frac{1}{2}+\frac{y}{2}-x}\Dfs(z)dzdy,
\end{equation*}
and we see that for this case we have that
\begin{equation*}
\p[(A\leqslant x)\cap((B-S)^++A> x)]=2n(n-1)\int_{(2x-1)^+}^{x}y^{n-2}\int_0^{\frac{1}{2}+\frac{y}{2}-x}\Dfs(z)dzdy.
\end{equation*}
\end{remark}

\begin{remark}
We see that as the size of an order, i.e.\ $n$ grows, the convergence to \Dfs\ is faster, and thus closer to the upper limit $F_1(x)=nx^{n-1}-(n-1)x^n$ (which was used as the first step of the iterations in Figure~\ref{fig4}. This is due to the fact that for $n\geqslant 2$,
\begin{equation*}
\left(\frac{n}{n+1}\right)^n\leqslant\left(\frac{n-1}{n}\right)^{n-1}.
\end{equation*}
Also observe that as the order size grows we have for the contraction constant that
$$
\lim_{n \to \infty} \left(\frac{n-1}{n}\right)^{n-1}=\lim_{n \to \infty} \left(1-\frac{1}{n}\right)^{n}=\frac{1}{e}.
$$
\end{remark}

\subsection{Shortest single direction after the nearest item}\label{ss:NI+short}
This strategy describes a slightly more sophisticated scenario. Rather than continue rotating in the same direction as the one where the nearest item lied, we first rotate to the nearest item and then choose the shortest single direction to pick the rest of the items, thus changing direction if needed.

There are several reasons why we study this strategy. First, the main point of our work, i.e.\ the analytical properties that emerge with each strategy, is clearly displayed again. Second, this strategy seems to be more efficient than the one previously examined, as it performs the same first step but then chooses a path which is at most as long as the previous case. The third reason is that this strategy is identical to the ``biggest-gap'' strategy (which we discuss in the following section) in case $n=1$ or $n=2$. The ``biggest-gap'' strategy reduces the travel time when picking only one order at a single carousel. In case one picks two orders at a pair of carousels, this is neither proven nor expected. It seems, however, interesting to study two strategies that overlap in the limiting cases (of the smallest order sizes possible).

In the results below we focus on $n \geqslant 3$, while we remark upon the remaining two cases (i.e.\ orders of size 1 or 2) at the end of the section. Our first task is to compute \Dfa, \Dfb\ and the conditional probability $\p[B\leqslant x | A=y]$. These are given in the following lemma.
\begin{lemma}\label{lem:condProb}
For $n\geqslant 3$, we have that:
\begin{equation}\label{eq:densB}
\Dfb(x)=1-(1-2x)^n,\quad 0\leqslant x \leqslant \frac{1}{2}.
\end{equation}
For $0\leqslant y \leqslant 1$ we have that:
\begin{equation}\label{eqA}
\Dfa(y)=2y^n+ny^{n-1}(1-y)-((2y-1)^+)^n-n((2y-1)^+)^{n-1}(1-y)
\end{equation}
and that
\begin{equation}\label{eq:condB2}
\p[B\leqslant x|A=y]=\begin{cases}
{\displaystyle \frac{2n(n-1)y^{n-2}x-n((y-2x)^+)^{n-1}+ny^{n-1}-4n(n-1)((2y-1)^+)^{n-2}x}{\dfa(y)}}, \\ \hfill \mbox{for } 0\leqslant x \leqslant \min\{\frac{1-y}{2},\frac{1}{4}\},\\
{\displaystyle \frac{1}{\dfa(y)}{\left(2n(n-1)y^{n-2}x+ny^{n-1}\right)}}, \hfill \mbox{for } \frac{1}{4}\leqslant x\leqslant\frac{1-y}{2},
\end{cases}
\end{equation}
where the right-hand side is presumed to be zero if the corresponding interval is invalid.
\end{lemma}
\begin{proof}
See Appendix~\ref{app:3}.
\end{proof}

We can now derive the integral equation that the sojourn time of the picker satisfies.
\begin{lemma}\label{lem:soj}
For $n\geqslant 3$ and $0 \leqslant x \leqslant 1$, the sojourn time of the picker for this strategy is given by the integral equation
\begin{align}
\nonumber
\Dfs(x)&=\Dfa(x)-\int_{(2x-1)^+}^x\int_0^{\frac{1+y}{2}-x}2n(n-1)y^{n-2}\Dfs(z)dzdy\\
\nonumber
       &\quad-\int_{\max\{(x-\frac{1}{4})^+,2x-1\}}^{x}\int_0^{\min\{\frac{1}{4}+y-x,\frac{1+y}{2}-x\}}\left(2n(n-1) ((3y-2x-2z)^+)^{n-2}\right)\Dfs(z)dzdy\\
&\quad+\int_{\max\{\frac{1}{2},2x-1\}}^{\max\{\frac{1}{2},x\}}\int_0^{\frac{1+y}{2}-x}\left(4n(n-1)(2y-1)^{n-2}\right) \Dfs(z) dzdy.
\end{align}
\end{lemma}
\begin{proof}
See Appendix~\ref{app:4}.
\end{proof}

Consider the complete metric space $\mathcal{L}$ with the sup norm. In Appendix~\ref{app:5} we prove the following theorem.
\begin{theorem}\label{th:3}
Consider \Dfa\ in \eqref{eqA}. For all $n \geqslant 3$ and $0 \leqslant x \leqslant 1$, the operator $\Omega:\mathcal{L}\rightarrow\mathcal{L}$ that is defined by
\begin{align*}
\Omega F(x)&=\Dfa(x)-\int_{(2x-1)^+}^x\int_0^{\frac{1+y}{2}-x}2n(n-1)y^{n-2}F(z)dzdy\\
       &\quad-\int_{\max\{(x-\frac{1}{4})^+,2x-1\}}^{x}\int_0^{\min\{\frac{1}{4}+y-x,\frac{1+y}{2}-x\}}\left(2n(n-1) ((3y-2x-2z)^+)^{n-2}\right)F(z)dzdy\\
&\quad+\int_{\max\{\frac{1}{2},2x-1\}}^{\max\{\frac{1}{2},x\}}\int_0^{\frac{1+y}{2}-x}\left(4n(n-1)(2y-1)^{n-2}\right) F(z) dzdy
\end{align*}
is a contraction mapping on $\mathcal{L}$ with contraction constant $C_n$ that is bounded from above by $\frac{11}{12}$.
\end{theorem}

\begin{remark}
Lemma~\ref{lem:soj}, and thus Theorem~\ref{th:3}, also hold in case $n=2$, provided that we account for indeterminate forms. Namely, one needs to be careful when plugging in $n=2$ in expressions such as $((3y-2x-2z)^+)^{n-2}$. In order to derive the correct expression, contrary to the usual convention, this quantity should be taken to be equal to 0, rather than 1, if the resulting term would be the indeterminate form $0^0$. Last, in the special case that orders consist of exactly one item, it can be shown that the sojourn time distribution is given by
$$
  \Dfs(x)=\begin{cases}
  \sin(2x) +\frac{1-\sin(1)}{\cos(1)}\cos(2x), &0\leqslant x\leqslant \frac{1}{2}\\1,&x\geqslant\frac{1}{2}.
\end{cases}
$$
\end{remark}

\subsection{Avoid the biggest gap}\label{ss:biggest_gap}

In this section we present some simulation results that compare the strategy in Section~\ref{ss:NI+short} to the ``biggest-gap'' strategy discussed earlier. One can mimic the derivations of the previous sections to show that the sojourn time satisfies a contraction mapping. Note that in case $n>2$, the sojourn time of the picker is not bounded by 1. The analysis will follow through with minor modifications; see e.g.\ Remark~\ref{rem:extendSpace}.

As noted before, for the special case where an order consists of two items, the ``avoid-the-largest-gap'' strategy coincides with the strategy where the carousel first rotates to the closest item and then chooses the shortest direction for the rest of the items. One can readily see that the sojourn time is bounded by 3/4 (which is a tight bound since we can achieve this sojourn time by having the starting position at 0, and the two items at 0.25 and 0.75). However, the analysis of this strategy is rather complex, and thus, we have chosen to depict some surprising results through simulations.

Evidently, for a single carousel, avoiding the biggest gap leads to the shortest travelling time for a single order. Thus, if one could pick the starting position, that would be at the end of the largest gap. However, picking multiple orders changes the dynamics; repositioning the carousel may lead to long travel times. Intuitively, it is not clear that this strategy is optimal for more carousels, where at least two orders are picked in a sequence, due to the dependencies between the preparation time of a carousel (i.e.\ the rotation time until the chosen starting position) and the picking time of an order. If we position the carousel ``ideally'', as would be the case for a single carousel, we minimise the time we need to pick the order at that carousel. However, that also reduces the available preparation time for the other machine, thus increasing the probability that upon return the picker will have to wait there until the other machine reaches the ``ideal'' starting point. On the other hand, if ``on average'' one can easily reach an endpoint of the largest gap, then this strategy should lead to the maximum throughput.

\begin{figure}
\includegraphics[width=0.49\textwidth]{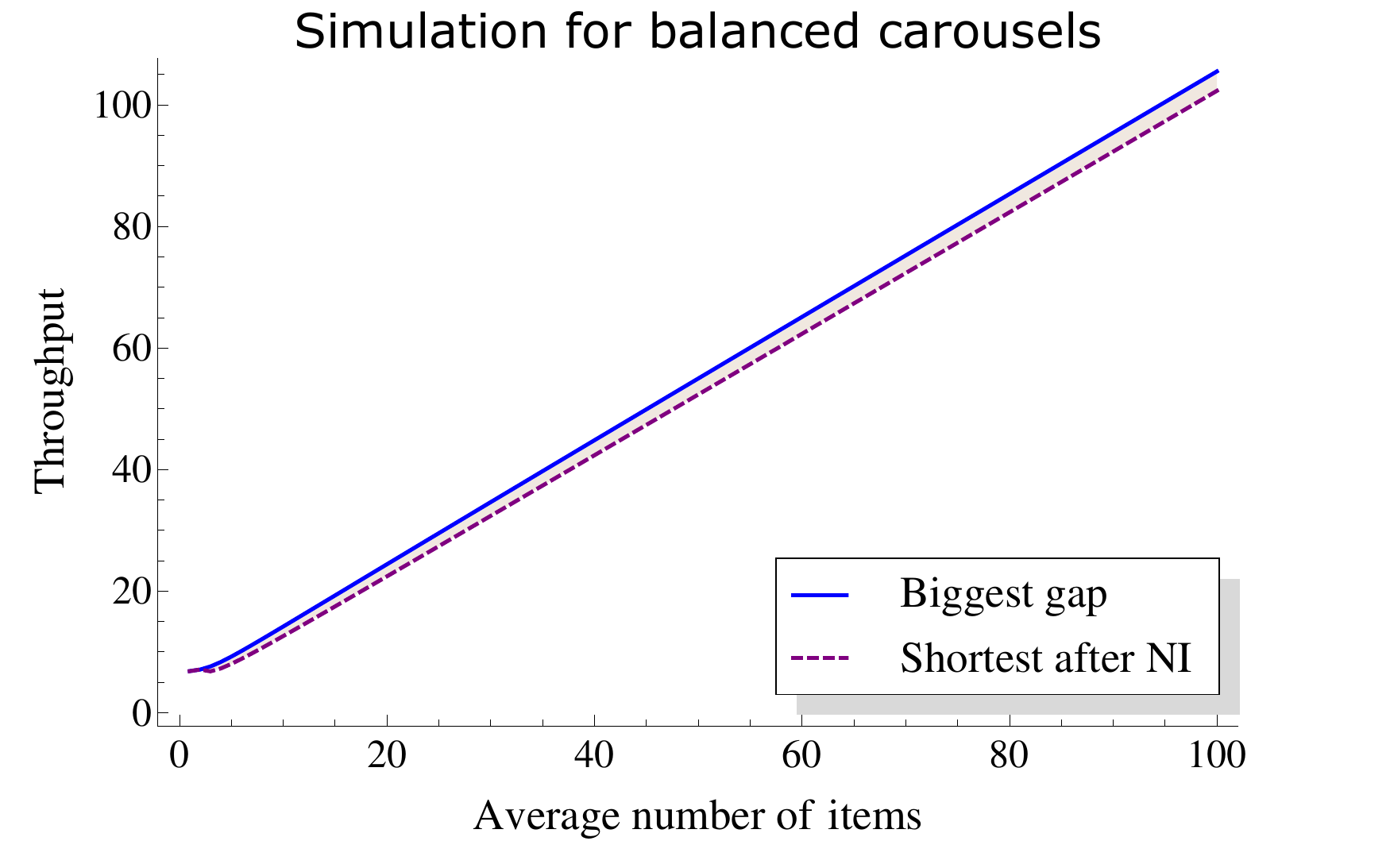}\hfill
\includegraphics[width=0.49\textwidth]{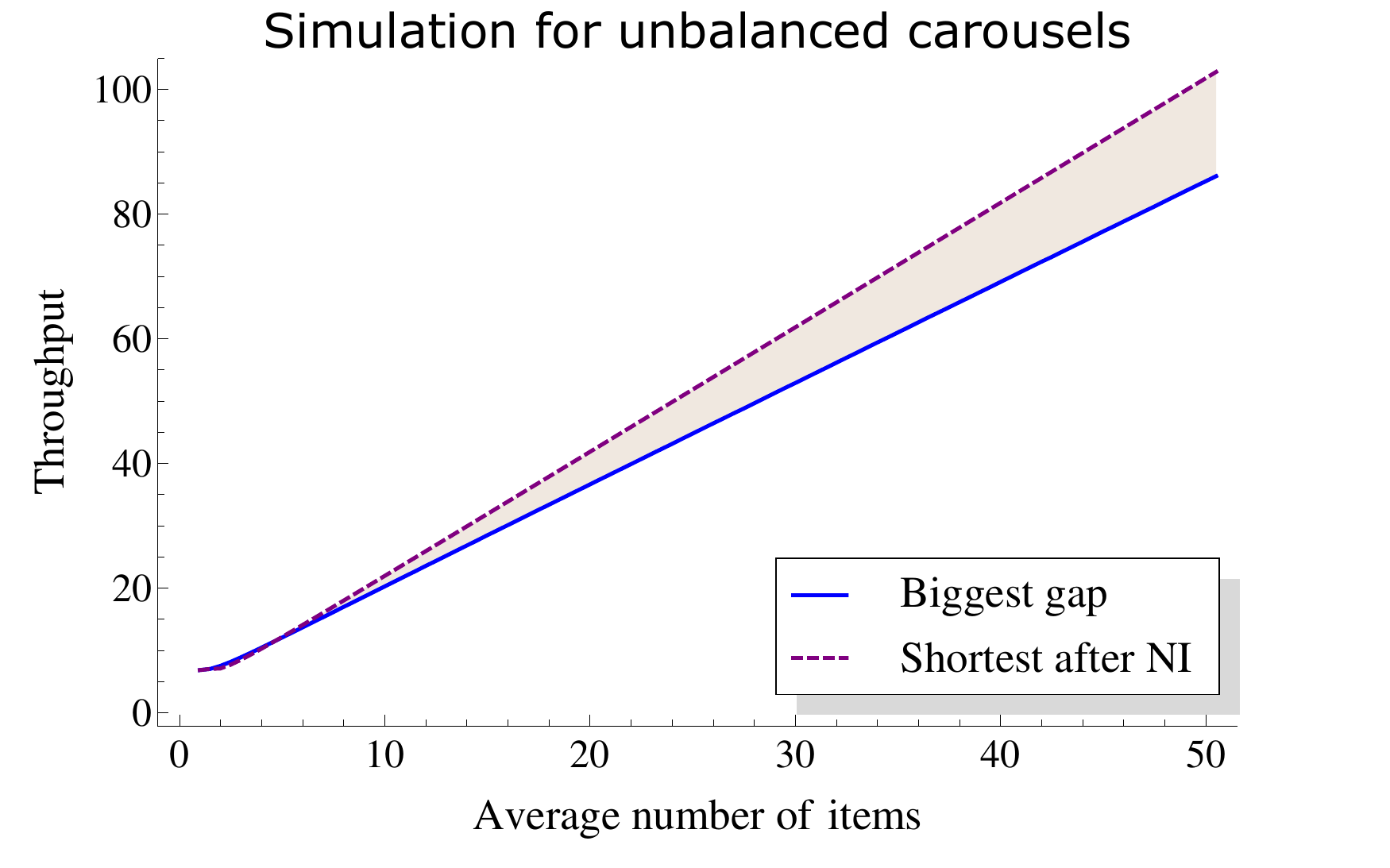}
\caption{Simulation results for balanced and unbalanced order sizes.\label{fixed}}
\end{figure}

In Figure~\ref{fixed} we plot the throughput of a pair of carousels for two scenarios. In both cases, each order has a fixed number of items. However, in the first case, the carousels are balanced: both carousels have orders of equal size $n$. In the second case the carousels are unbalanced: all orders on the first carousel consist of 1 item, while the orders on the second carousel consist of $n$ items. The simulation is done for $n=1,2,\ldots,100$ for both strategies we consider. The horizontal axis gives the average size of a random order.

As seen from Figure~\ref{fixed}, avoiding the biggest gap leads to a bigger throughput if the carousels are balanced. Avoiding the biggest gap is still optimal in the unbalanced case, provided that the difference between the order sizes in the two carousels is limited. Already from the moment the second carousel has order sizes of 9 items or more, the strategy in Section~\ref{ss:NI+short} outperforms the ``biggest-gap'' strategy. It is intuitively clear that for a large order size on the second carousel, avoiding the biggest gap will not be optimal. The biggest gap can lie relatively far from the origin, and thus, as expected, moving to the nearest item and choosing the shortest direction from that point leads to a bigger throughput due to the smaller preparation time. What is not intuitively clear is why the difference in the preparation time becomes significant only from the point that orders on the second carousel are 9 times as large as orders on the first one.

\begin{figure}
\begin{minipage}{\textwidth}
\includegraphics[width=0.49\textwidth]{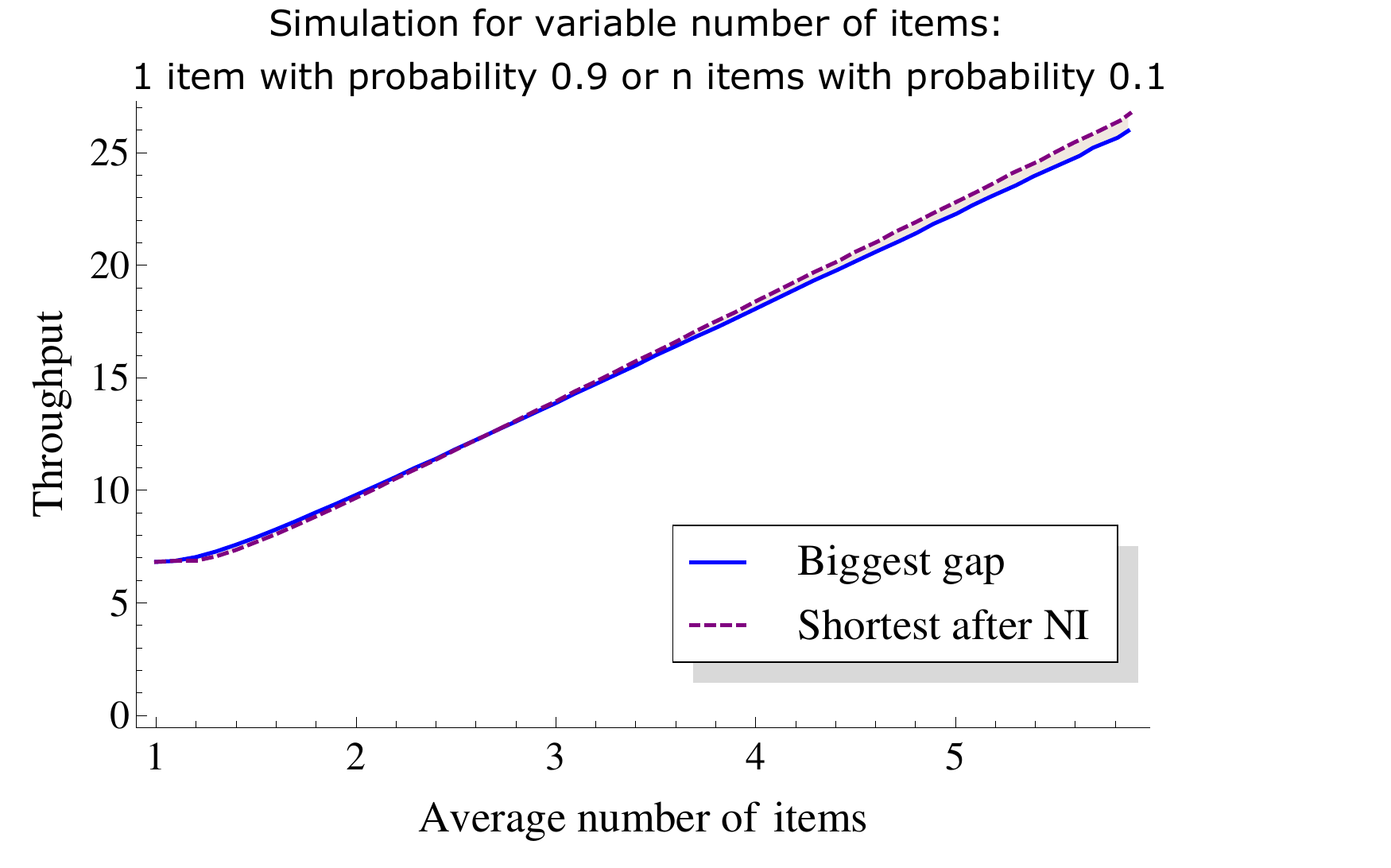}\hfill
\includegraphics[width=0.49\textwidth]{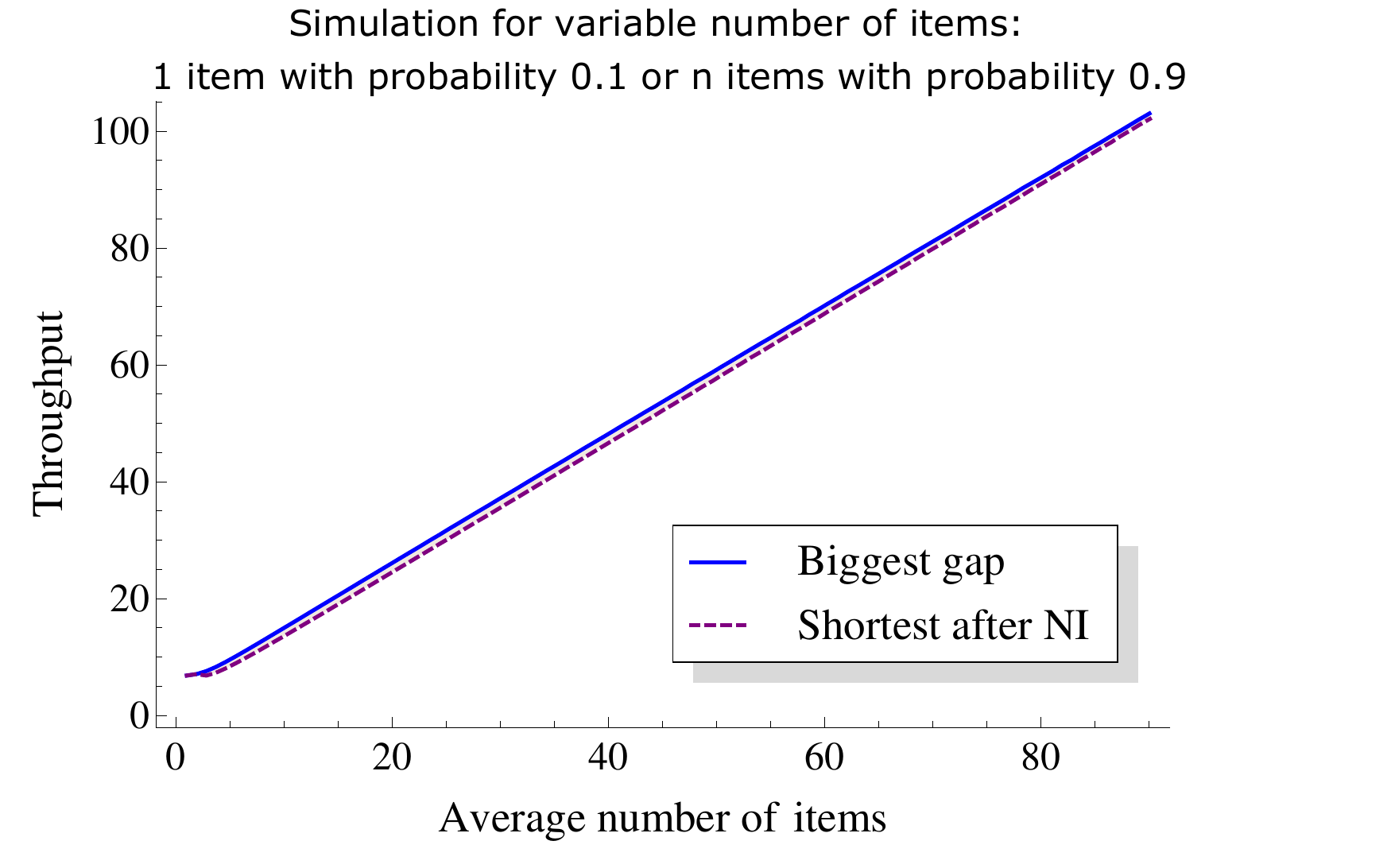}
\end{minipage}
\begin{minipage}{\textwidth}
\begin{center}
\includegraphics[width=0.49\textwidth]{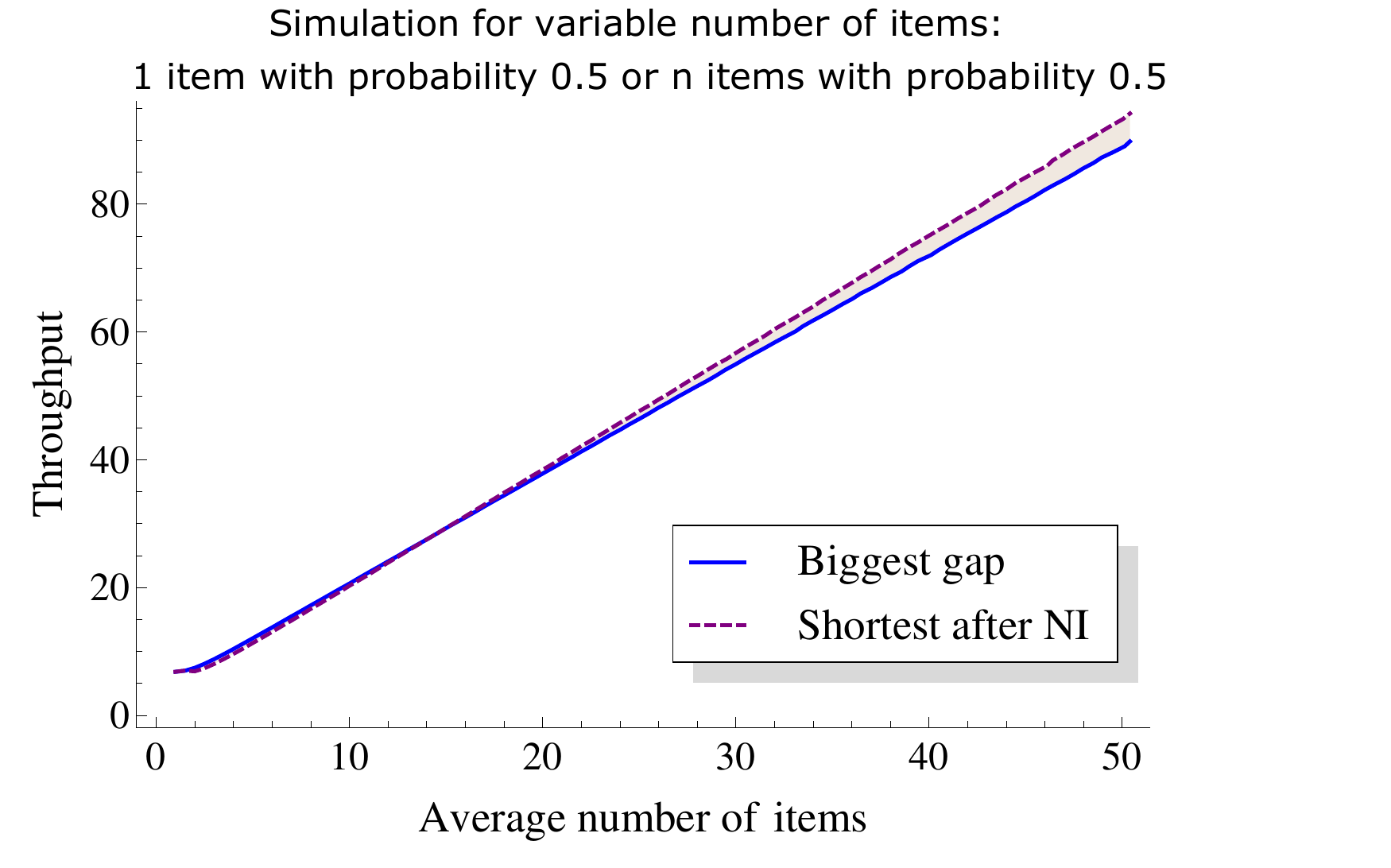}
\end{center}
\end{minipage}
\caption{Simulation results for variable order sizes. \label{variable}}
\end{figure}

Similarly to Section~\ref{s:variable}, we now investigate other sources of variability. Namely, we now allow for variable order sizes. Each order, on both carousels, consists of 1 item with probability $p$ and of $n$ items with probability $1-p$. The simulation is done for $n=1,2,\ldots,100$ and $p\in\{0.1,0.5,0.9\}$, for both the ``biggest-gap'' strategy and the ``shortest-direction-after-the-nearest-item'' strategy. Results are given in Figure~\ref{variable}. Apparently, the ``biggest-gap'' strategy performs relatively worse for a larger fraction of small orders, which is somewhat surprising.

Although comparing fixed order sizes with variable order sizes may be impossible, the conclusion from Figure~\ref{fixed} was that when the variability between the two carousels was low, the ``biggest-gap'' strategy outperformed the ``shortest-direction-after-the-nearest-item'' strategy. However, as seen in the case examined in the top-right part of Figure~\ref{variable}, the ``biggest-gap'' strategy can outperform the ``shortest-direction-after-the-nearest-item'' strategy even when variability is high; admittedly we now talk about variability in the size of a single order rather than the order size between the two carousels.

What is evident from Figures~\ref{fixed} and~\ref{variable} is that predicting the optimal strategy for a pair of carousels is not very intuitive. However, time-consuming simulations can be avoided by simply iterating the integral equation \eqref{eq2}, as the sojourn time distribution can be shown to satisfy a contraction mapping in $\mathcal{L}$ under a variety of strategies, such as the ones considered in this paper.

\section{Simulation results}\label{s:sim}

In the previous sections we have mainly focused on showing the contraction property for three out of the four strategies presented, while we claimed that a similar result can be derived for the last strategy (the ``biggest-gap'' strategy), at the cost of lengthier derivations. An indication about the validity of our claim was given by the fact that the ``biggest-gap'' strategy coincides with the one described in Section~\ref{ss:NI+short} for orders of size 2.

In this section we give some simulation results in an effort to further compare the four strategies presented so far, together with three additional reasonable strategies. As an optimal strategy for two carousels is not yet known, the results presented in this section offer some first indication on how each strategy that we have considered performs. We do not present exhaustive simulation results. We rather focus on a few main conclusions.

The strategies we consider in this section are as follows:
\begin{enumerate}
\item \label{str:OD_NI}
Unidirectional carousels: the preparation phase consists of rotating to the first item in the direction of rotation (Section~\ref{s:one_dir}).
\item Unidirectional carousels: the preparation phase consists of rotating to the item \textit{after} the biggest gap. Thus, when picking the rest of the items, the biggest gap will be avoided. This is a sensible strategy if the probability that the preparation can be completed during the pick time of the previous order is relatively large.
\item \label{str:TD_NI}
Bidirectional carousels: the preparation phase consists of rotating to the nearest item -- whether that lies at the right or at the left of the origin -- and later on continuing in that direction in order to collect the remaining items. The carousel, thus, does not reverse its direction during the picking of one order (Section~\ref{s:NI}).
\item \label{str:NI+short}
Bidirectional carousels: the preparation phase consists of rotating to the nearest item; from that position, the remaining items will be collected along the shortest path to the right or the left. Thus, the carousel may reverse its direction once after collecting the first item (Section~\ref{ss:NI+short}).
\item \label{str:TW_BG}
Bidirectional carousels: the carousel rotates to the nearest item from the two that lie on the right and the left of the biggest gap. The rest of the items will be collected in the direction that avoids traversing the biggest gap, which minimises the travel time per order (Section~\ref{ss:biggest_gap}).
\item \label{str:2nd_best}
Bidirectional carousels: the carousel rotates to the item that has the property that the second item in that direction is closest to the origin; the rest of the items are collected along the shortest path. This strategy is motivated by the fact that the optimal policy admits at most one turn \cite{bartholdi86}, and for small order sizes, the second item to the origin offers an indication to the direction in which most items lie.
\item \label{str:TW_maybe_BG}
Bidirectional carousels: rotate to the nearest endpoint of the biggest gap, provided that this can be completed during the pick time of the previous order; then collect the remaining items by avoiding the biggest gap. Should it be impossible to reach an endpoint of the biggest gap within time, check if reaching an endpoint of the second biggest gap within time is possible, and if so move there. If not, check if the carousel can reach an endpoint of the third biggest gap and so on. If no item can be reached within the pick time of the previous order, simply remain at the origin and choose the shortest direction (either the one to the right or the one to the left).
\end{enumerate}

\begin{figure}
\includegraphics[width=0.49\textwidth]{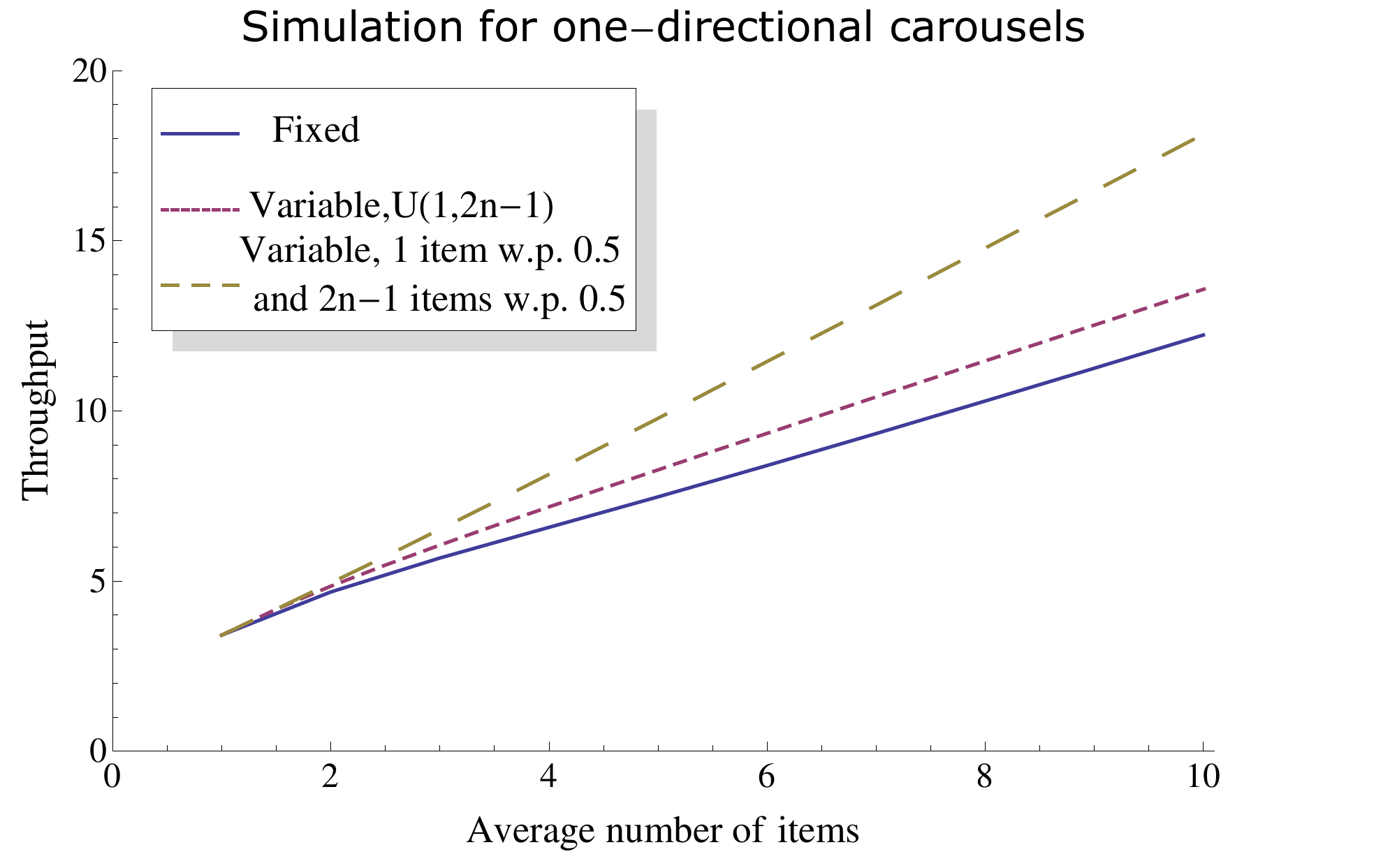} \hfill
\includegraphics[width=0.49\textwidth]{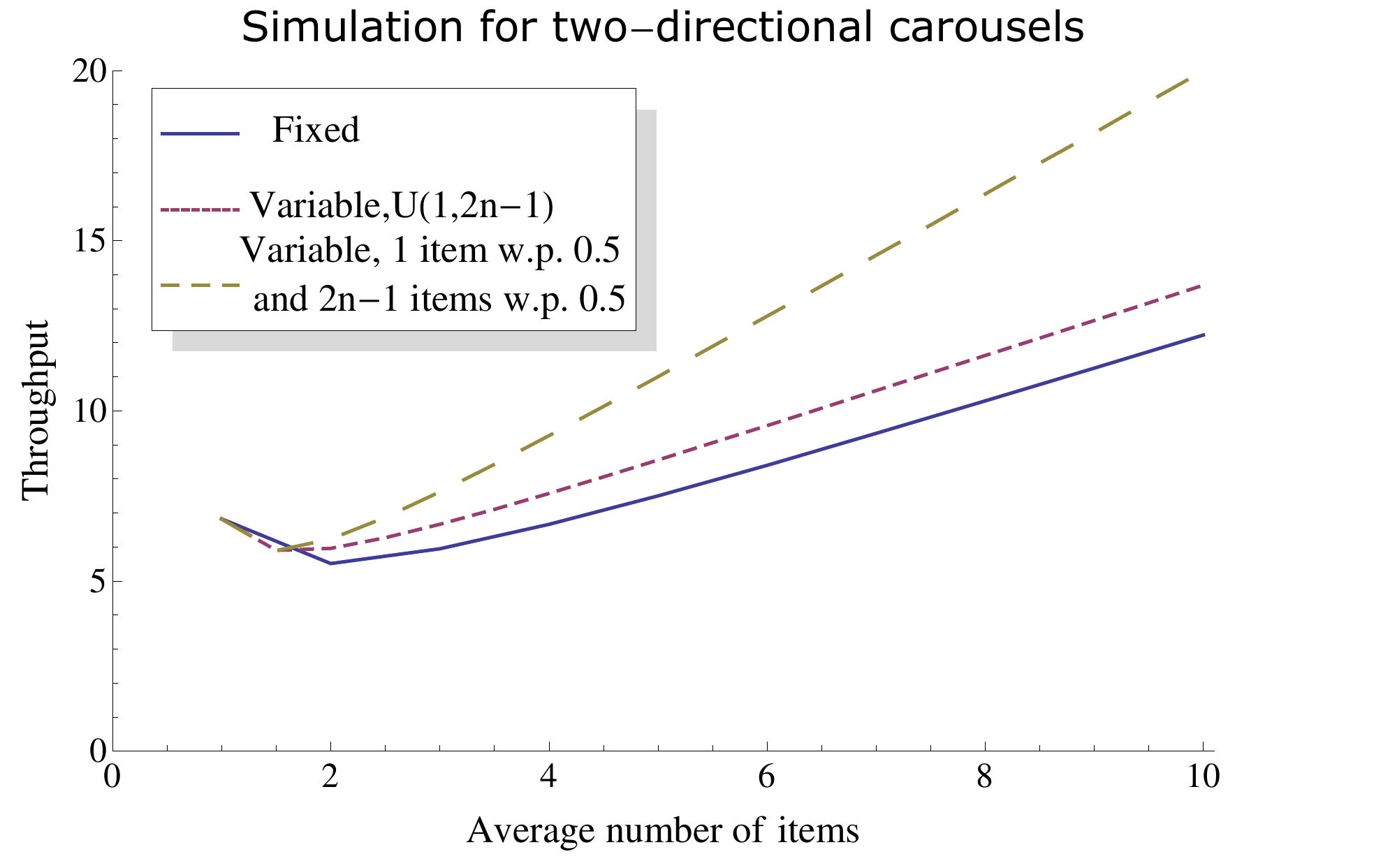}
\caption{Comparing the throughput of strategies \ref{str:OD_NI} and~\ref{str:TD_NI} for three order size scenarios. The order size is either fixed, or uniformly distributed on the integers between 1 and $2n-1$, or either 1 or $2n-1$ with probability 0.5 each. \label{fig:NI_comparison}}
\end{figure}

In Figure~\ref{fig:NI_comparison}, we examine the gain in throughput for a pair of unidirectional carousels vs.\ a pair of bidirectional carousels, both for fixed and for variable order sizes. Namely, we compare strategies~\ref{str:OD_NI} and~\ref{str:TD_NI}, which are very similar to one another. Their main difference is that the rotation time to the first item for strategy~\ref{str:TD_NI} is on average half of that of strategy~\ref{str:OD_NI}. It is thus expected that for $n=1$ the sojourn time of the second strategy is half of that of the first one, and that as $n$ increases, the difference between the two strategies decreases, as the benefit from a shorter preparation time becomes smaller, since the distance to the first item becomes shorter.

Another conclusion that may be drawn from this picture is that variability in order sizes leads to bigger throughput. We examine three scenarios: fixed order sizes, uniformly distributed order sizes varying from 1 item to $2n-1$ items, and order sizes that are either precisely equal to 1 with probability 1/2 or equal to $2n-1$ (thus giving an average order size of $n$ for all three scenarios). The variance of the first distribution is zero; the second has variance equal to $\frac{1}{3} n (n-1)$, and the last has variance equal to $(n-1)^2$, which is bigger than the other two distributions for all (integer) order sizes $n$. We see that not only is the throughput higher as the variance increases, but that also the the rate of increase is higher for higher variances.

\begin{figure}
\begin{center}
\includegraphics[width=0.9\textwidth]{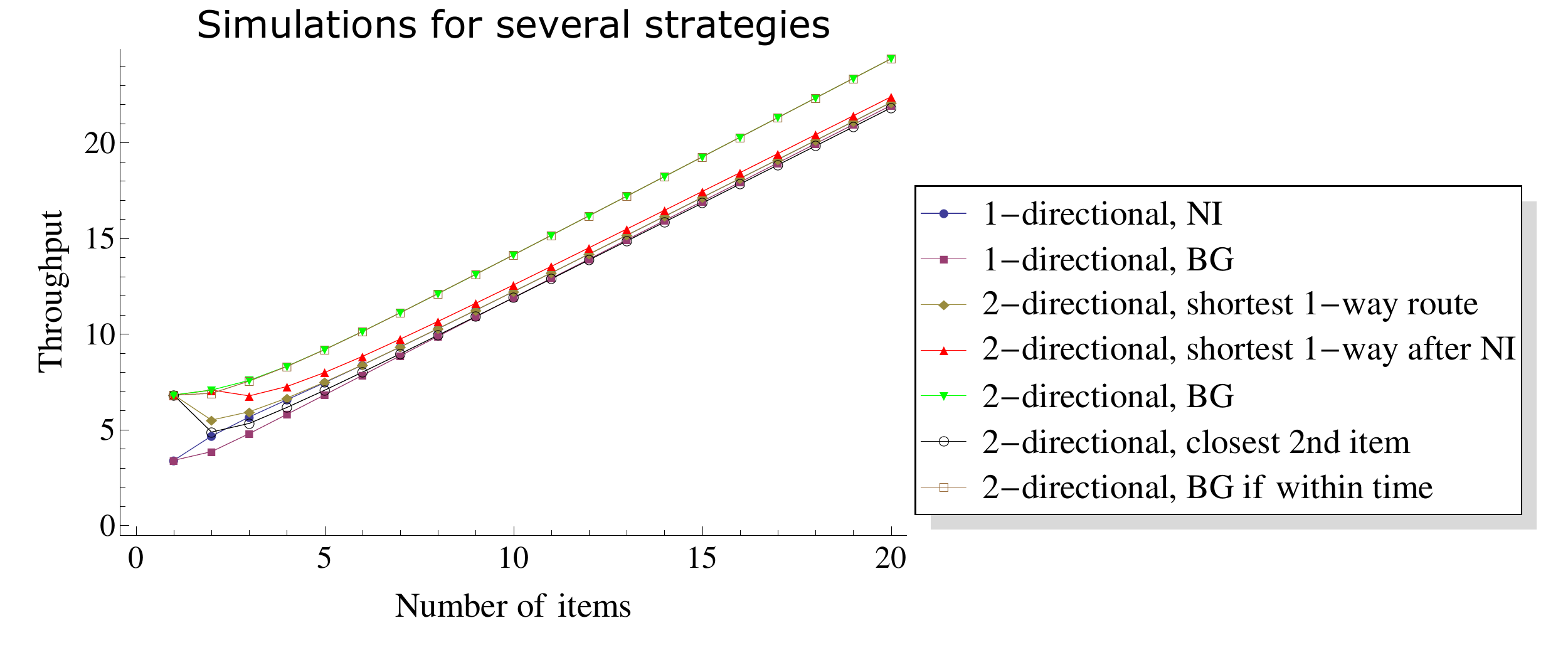}
\end{center}
\caption{Comparing the throughput of all strategies for fixed order sizes.\label{fig:comp}}
\end{figure}

In Figure~\ref{fig:comp}, we compare the performance of all seven strategies for a fixed order size. We see that strategies~\ref{str:TW_BG} and~\ref{str:TW_maybe_BG} perform best in this case. The reason is that the travel time is minimised since we avoid traversing the biggest gap. Keep in mind though that this conclusion strongly depends on the fact that we have balanced order sizes. As we saw in Figures~\ref{fixed} and~\ref{variable}, this conclusion does not hold in general. Strategy~\ref{str:TW_BG} is not uniformly better than strategy~\ref{str:NI+short}. In case of balanced order sizes, we see that the benefit of strategies~\ref{str:TW_BG} and~\ref{str:TW_maybe_BG} increases as the order size increases; as $n$ increases, the probability that the picker will be busy long enough with the previous order to reach an endpoint of the biggest gap in time increases.

Strategy \ref{str:2nd_best} looks at the distance to the second item from the origin in both directions. Although one can easily imagine situations where this strategy would perform the best, it is by far the least useful strategy from all strategies with motion in two directions. Moreover, as the (fixed) order size increases, this strategy performs worse than all others.

As mentioned earlier, for $n=1$ and $n=2$ strategies \ref{str:NI+short} and \ref{str:TW_BG} are identical; cf.\ Figure~\ref{fig:comp}. For fixed order sizes, strategy \ref{str:NI+short} performs worse than strategies \ref{str:TW_BG} and \ref{str:TW_maybe_BG}, although better than the rest of the picking strategies we have considered. A reason for this is that the preparation time $B$ is minimised (as is the case for strategy \ref{str:TD_NI}, only that now once chooses the shortest path after the first item, rather than continuing in the same direction, which naturally increases the throughput). However, as at least $n-1$, if not all, items are collected in one direction, one may have a long pick time. Consider for example the situation where 3 items are positioned on points 0.1, 0.2, and 0.8, where the origin is located at 0 and directly across the origin we have point 0.5. In this situation, strategy \ref{str:NI+short} would first lead to the item on 0.1, and then from that point on, the carousel would traverse the interval $[0.1,0.8]$ (travel time of 0.7), which is shorter than going to point 0.2 in the other direction (travel time of 0.9). In this scenario, strategies \ref{str:TW_BG}, \ref{str:TW_maybe_BG}, and even \ref{str:2nd_best} would perform better.

An interesting conclusion is drawn when studying strategies \ref{str:NI+short}, \ref{str:TW_BG}, and \ref{str:TW_maybe_BG}. The latter is a compromise between the first two. The idea is that when trying to avoid the biggest gap by strategy \ref{str:TW_BG}, one may have a too long preparation time, while when trying to minimise the preparation time, as in strategy \ref{str:NI+short}, one may end up with a too long pick time. Strategy \ref{str:TW_maybe_BG} uses more information (i.e.\ the pick time of the previous order) to decide whether one can reach an end point of the biggest gap in time. As is evident from Figure~\ref{fig:comp}, strategies \ref{str:TW_BG} and \ref{str:TW_maybe_BG} are almost identical and have the highest throughput among all strategies we studied for fixed order sizes, with strategy \ref{str:TW_BG} performing better on average. The added information of strategy \ref{str:TW_maybe_BG} offers little improvement (about 0.04\% more than strategy \ref{str:TW_BG} in the best case). Moreover, for small $n$, avoiding the biggest gap at all times performs better. Intuitively, one can see the merit of having the picker wait for the remaining preparation time in order to reach an end point of the biggest gap, in particular since if the order size is small, then the biggest gap is expected to be substantial with respect to a remaining preparation time.

Another interesting observation is that avoiding the biggest gap, which leads to the highest throughput for fixed order sizes and bidirectional carousels, performs worse than going to the nearest item for unidirectional carousels. Combining this information, with the results in Figure~\ref{variable}, we conclude that the performance of a strategy is strongly influenced by the variability of the order size on each carousel, the variability of (average) order sizes between the two carousels, and whether a carousel can rotate in two directions or not. No strategy among the ones we have considered outperforms the other ones with respect to all three considerations, although avoiding the biggest gap or choosing the shortest direction after moving to the nearest item are two strategies that are usually reasonable.

\phantomsection
\pdfbookmark[1]{Acknowledgement}{ack}
\section*{Acknowledgement}
The authors would like to thank Jos Brands for initiating the present study and pointing us to the fixed point equation for the sojourn time
distribution.

\phantomsection
\pdfbookmark[1]{References}{refs}
\bibliographystyle{apt}
\bibliography{maria}

\newpage
\section*{Appendix}
\appendix

\section{Proof of Lemma~\ref{lem:1}} \label{app:1}
Observe that when $y\geqslant 2x-1$ or when $0\leqslant x\leqslant\frac{1}{2}$ and $0\leqslant y\leqslant x$ we have that $\frac{1+y}{2}-x \geqslant 0$. From \eqref{eq2}, we have for $0 \leqslant x \leqslant 1$:
\begin{align*}\Dfs(x)&=\int_{0}^{x}\int_{0}^{1}\p[B\leqslant x-y+z|A=y]\dfa(y)\dfs(z)dzdy\\
&=\int_{0}^{(2x-1)^+} \int_{0}^{1} \dfa(y) \dfs(z) dzdy+\int_{(2x-1)^+}^{x}\int_{0}^{a}\p[B\leqslant x-y+z|A=y] \dfa(y) \dfs(z) dzdy \\
&\quad+\int_{(2x-1)^+}^{x}\int_{a}^{1}\dfa(y)\dfs(z)dzdy.
\end{align*}

This is true for the following reasons. For $x\leqslant \frac12$, the first integral is equal to zero. For $x>\frac12$, we have that $x-y+z>\frac{1-y}{2}$ and thus the conditional probability we had is equal to one. A similar reasoning can be given for the last integral. Thus, we examine two separate cases, depending on the range of $x$:

For $0\leqslant x \leqslant \frac{1}{2}$, choose $a=\frac{1+y}{2}-x$ and observe that $a\geqslant 0$  and $2x-1\leqslant 0$. We now have from \eqref{eq3} and \eqref{eq:condB},
\begin{align*}
\Dfs(x)&=n(n-1)\left[2\int_{0}^{x}\int_{0}^{a}(x-y+z)y^{n-2}\dfs(z)dzdy +\int_{0}^{x}\int_{a}^{1}\left(y^{n-2}-y^{n-1}\right)\dfs(z)dzdy\right]\\
&=n(n-1)\left[2\int_{0}^{x}y^{n-2}\left((x-y)\Dfs(a)+\int_{0}^{a}z\dfs(z)dz\right)dy+\int_{0}^{x}(y^{n-2}-y^{n-1}) (1-\Dfs(a))dy\right]\\
&=nx^{n-1}-(n-1)x^n+n(n-1)\left[2\int_{0}^{x}xy^{n-2}\Dfs(a)dy-\int_{0}^{x}y^{n-1}\Dfs(a)dy\right]\\
&\quad+n(n-1)\left[2\int_{0}^{x}y^{n-2}\left(a\Dfs(a)-\int_{0}^{a}\Dfs(z)dz\right)dy-\int_{0}^{x}y^{n-2}\Dfs(a)dy \right]\\
&=nx^{n-1}-(n-1)x^n\\
&\quad+n(n-1)\left[\int_{0}^{x}y^{n-2}\Dfs(a)(2x-y+2a-1)dy-2\int_{0}^{x}y^{n-2}\int_{0}^{a}\Dfs(z)dzdy\right]
\end{align*}
Now, for $a=\frac{1+y}{2}-x$ we have that the first integral is equal to zero, and thus,
\begin{equation}\label{eq42}
\Dfs(x)=nx^{n-1}-(n-1)x^n-2n(n-1)\int_{0}^{x}y^{n-2}\int_{0}^{\frac{1}{2}+\frac{y}{2}-x}\Dfs(z)dzdy.
\end{equation}

For $\frac{1}{2}\leqslant x\leqslant 1$ we have that $2x-1\leqslant x$ and that
\begin{multline*}
\Dfs(x)=n(n-1)\bigg[\int_{0}^{2x-1}\int_{0}^{1}(1-y)y^{n-2}\dfs(z)dzdy+2\int_{2x-1}^{x}\int_{0}^{a}(x-y+z)y^{n-2} \dfs(z)dzdy\\
+\int_{2x-1}^{x}\int_{a}^{1}(y^{n-2}-y^{n-1})\dfs(z)dzdy\bigg].
\end{multline*}
Simple integration now yields
\begin{multline*}
\Dfs(x)=nx^{n-1}-(n-1)x^n+n(n-1)\bigg(\int_{2x-1}^{x}y^{n-2}\Dfs(a)(2x-y+2a-1)dy  \\-2\int_{2x-1}^{x}y^{n-2}\int_{0}^{a}\Dfs(z)dzdy\bigg),
\end{multline*}
which again for the particular choice of $a$ we have made is equal to
\begin{equation}\label{eq5}
\Dfs=nx^{n-1}-(n-1)x^n-2n(n-1)\int_{2x-1}^{x}y^{n-2}\int_{0}^{\frac{1}{2}+\frac{y}{2}-x}F(z)dzdy.
\end{equation}
Thus, for $0\leqslant x\leqslant 1$, we have \eqref{eq:int} by combining \eqref{eq42} and \eqref{eq5}, which is what we wished to prove. As an aside remark, after some computations, one can show that
\begin{equation*}
\Dfs['](x)=n(n-1)x^{n-2}-n(n-1)x^{n-1}-x^{n-2}\int_0^{\frac{1-x}{2}}\Dfs(z)dz+ \int_{(2x-1)^+}^{x}y^{n-2}\Dfs\left(\tfrac{1+y}{2}-x\right)dy.
\end{equation*}

\section{Proof of Theorem~\ref{th:2}} \label{app:2}
Let $F_1(x)$ and $F_2(x)$ be continuous functions in $\mathcal{L}$ with $|F_1(x)-F_2(x)|\leqslant \delta$ for all $x\in [0,1]$. We then have that
\begin{align*}
|\Omega F_1(x)-\Omega F_2(x)|&\leqslant 2n(n-1)\int_{(2x-1)^+}^{x}y^{n-2}\int_0^{\frac{1}{2}+\frac{y}{2}-x}|F_1(z)-F_2(z)|dzdy\\
&\leqslant 2n(n-1)\int_{(2x-1)^+}^{x}y^{n-2}\int_0^{\frac{1}{2}+\frac{y}{2}-x}\delta dzdy\\
&= \delta
2n(n-1)\left(\frac{x^{n-1}}{2(n-1)}+\frac{x^n}{2n}-\frac{x^n}{n-1}-\frac{((2x-1)^+)^{n-1}}{2(n-1)}\right)\\
&\quad-\delta2n(n-1)\left(\frac{((2x-1)^+)^n}{2n}-\frac{x((2x-1)^+)^{n-1}}{n-1}\right)\\
&=\delta(x^{n-1}(n+nx-x-2nx)-((2x-1)^+)^{n-1}(n+(n-1)(2x-1)^+-2nx)).
\end{align*}
Define now
$c(x)=x^{n-1}(n+nx-x-2nx)-((2x-1)^+)^{n-1}(n+(n-1)(2x-1)^+-2nx)$ and observe that
\begin{align*}
c(x)&\leqslant x^{n-1}(n-x-nx)-((2x-1)^+)^{n-1}(n+2nx-n-2x+1-2nx)\\
&=x^{n-1}(n-x-nx)+((2x-1)^+)^{n}\\
&\leqslant x^{n-1}(n-x-nx)+x^n=nx^{n-1}(1-x).
\end{align*}

Observe now that $nx^{n-1}(1-x) \geqslant 0$ for all $x\in [0,1]$ and that the maximum value is given at the point $x=\frac{n-1}{n}$ and is equal to
\begin{equation*}
n\left(\frac{n-1}{n}\right)^{n-1}(1-\frac{n-1}{n})=\left(\frac{n-1}{n}\right)^{n-1}.
\end{equation*}
Moreover, we have that $0<\left(\frac{n-1}{n}\right)^{n-1}<1$ for all $n\geqslant 2$, and thus $\Omega$ is a contraction mapping on $\mathcal{L}$ with contraction constant $\left(\frac{n-1}{n}\right)^{n-1}$.

\section{Proof of Lemma~\ref{lem:condProb}}\label{app:3}
For the distribution of $B$, observe that $\Dfb(x)=1-\p[B>x]=1-(1-2x)^n$, since the items are uniformly distributed and since if $x$ is the distance to the nearest item (say at the right of the origin), then the distance to the first item from the other side (say the left) is at least $x$, and thus all $n$ items are located in an interval of length $1-2x$. As we first move to the nearest item, we have that $x$ cannot be more than $1/2$.

In order to derive \Dfa, we first need to derive the conditional probability $\p[A\leqslant y|B=x]$. We immediately have that for $\frac{1}{4}\leqslant x\leqslant\frac{1}{2}$ and $0\leqslant y\leqslant 1-2x$,
\begin{equation}\label{eq:a|b1}
\p[A\leqslant y|B=x]=\left(\frac{y}{1-2x}\right)^{n-1}.
\end{equation}

To see this, recall that the items are uniformly distributed on the carousel and that since we first go to the nearest item, the distance to both items to the left and to the right of the starting position is at least $x$. Moreover, if $\frac{1}{4}\leqslant x \leqslant \frac{1}{2}$, we have that the direction of rotation after the first item will remain the same.

\begin{figure}[h!]
\begin{center}
\includegraphics[width=5cm]{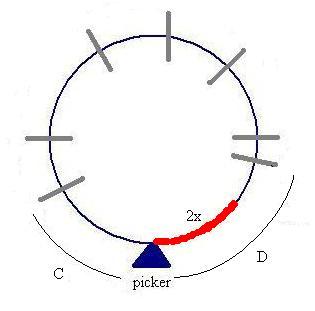}
\caption{If $C\leqslant D$, the remaining $n-1$ items will be picked by turning to the left.\label{fig7}}
\end{center}
\end{figure}

The computation for $0\leqslant x\leqslant \frac{1}{4}$ is a bit more complicated. The reason is that after reaching the nearest item, which lies at a distance $B=x$, the remaining $n-1$ items may have to be picked either by maintaining or by reversing the rotation direction. After reaching the nearest item, we know that in some direction (say, without loss of generality that this is the right one) there are no more items for at least a distance of $2x$. Let $C$ be the distance to the nearest item at the left, and $D$ be the distance to the nearest right item; cf.\ Figure~\ref{fig7}.

We now condition on the position on the nearest \emph{left} item. Thus, for $0\leqslant x\leqslant \frac{1}{4}$ and  $0\leqslant y \leqslant 1-2x$ we have that
\begin{align}\label{eq:a|b2}
\nonumber \p[A\leqslant y|B=x]&=\int_0^{2x}\p[A\leqslant y|B=x,C=z]f_{C|B=x}(z)dz\\
\nonumber &\quad+ \p[A\leqslant y|B=x,C\geqslant 2x]\p[C\geqslant 2x|B=x]  \\
\nonumber &=\int_0^{\min\{y,2x\}}\left(\frac{y-z}{1-2x-z}\right)^{n-2}\left(\frac{n-1}{1-2x}\right)\left(\frac{1-2x-z}{1-2x} \right)^{n-2}dz\\
\nonumber &\quad+\left(\frac{1-4x}{1-2x}\right)^{n-1}\left[2\left(\frac{(y-2x)^+}{1-4x}\right)^{n-1}-\left(\frac{(2y-1)^+}{1-4x} \right)^{n-1}\right]\\
\nonumber &=\int_0^{\min\{y,2x\}}(n-1)\frac{(y-z)^{n-2}}{(1-2x)^{n-1}}dz+2\left(\frac{(y-2x)^+}{1-2x}\right)^{n-1}- \left(\frac{(2y-1)^+}{1-2x}\right)^{n-1}\\
&=\left(\frac{y}{1-2x}\right)^{n-1}+\left(\frac{(y-2x)^+}{1-2x}\right)^{n-1}-\left(\frac{(2y-1)^+}{1-2x}\right)^{n-1}.
\end{align}
As a sanity check, we indeed verify that $\p[A\leqslant 0|B=x]=0$ and that $\p[A\leqslant 1-2x|B=x]=1$.

We proceed now with the distribution of $A$. For $0\leqslant y \leqslant \frac{1}{2}$ we have
\begin{align*}
&\Dfa(y)=\int_0^{\frac{1}{2}}\p[A\leqslant y|B=x]\dfb(x)dx\\
&=\int_0^{\frac{y}{2}}2n(y-2x)^{n-1}dx-\int_0^{\frac{y}{2}}2n((2y-1)^+)^{n-1}dx +\int_0^{\frac{1-y}{2}}2ny^{n-1}dx+\int_{\frac{1-y}{2}}^{\frac{1}{2}}2n(1-2x)^{n-1}dx\\
&=2y^n+ny^{n-1}(1-y).
\end{align*}
Furthermore, for $\frac{1}{2}\leqslant y\leqslant 1$ we have
\begin{align*}
&\Dfa(y)=\int_0^{\frac{1}{2}}\p[A\leqslant y|B=x]\dfb(x)dx\\
&=\int_0^{\frac{1-y}{2}}2n(y-2x)^{n-1}dx-\int_0^{\frac{1-y}{2}}2n(2y-1)^{n-1}dx+ \int_0^{\frac{1-y}{2}}2ny^{n-1}dx+\int_{\frac{1-y}{2}}^{\frac{1}{2}}2n(1-2x)^{n-1}dx\\
&=2y^n+ny^{n-1}(1-y)-(2y-1)^n-n(2y-1)^{n-1}(1-y).
\end{align*}
We see that \Dfa\ is continuous at $y=\frac{1}{2}$, and that aggregating these results leads to \eqref{eqA} for $0\leqslant y \leqslant 1$.

For $0\leqslant y \leqslant 1$, the density of the picking time is given by
\begin{align}\label{fay}
\dfa(y)&=2ny^{n-1}-ny^{n-1}+n(n-1)y^{n-2}(1-y)-2n((2y-1)^+)^{n-1}\nonumber\\
&\quad+n((2y-1)^+)^{n-1}-2n(n-1)((2y-1)^+)^{n-2}(1-y)\nonumber\\
&=ny^{n-2}(y+(n-1)(1-y))-n((2y-1)^+)^{n-2}(2(n-1)(1-y)+2y-1).
\end{align}

For the conditional distribution $\p[B\leqslant x | A=y]$ we proceed as follows. Denote by $f_{X|Y=y}(x)$ the conditional density of the random variable $X$ conditioned on the event $Y=y$. It is known that
\begin{equation}\label{eq:condDens}
f_{B|A=y}(x)=\frac{f_B(x)}{f_A(y)}f_{A|B=x}(y).
\end{equation}
Differentiating \eqref{eq:a|b1} and \eqref{eq:a|b2} we obtain (for $0\leqslant y \leqslant 1-2x$) that the conditional density $f_{A|B}$ is given by
\begin{equation}\label{eq:condA}
f_{A|B=x}(y) =\begin{cases}
\frac{n-1}{(1-2x)^{n-1}}\left(y^{n-2}+((y-2x)^+)^{n-2}-2((2y-1)^+)^{n-2}\right), &0\leqslant x\leqslant \frac{1}{4},\\
(n-1)\frac{y^{n-2}}{(1-2x)^{n-1}}, &\frac{1}{4}\leqslant x\leqslant \frac{1}{2}.
  \end{cases}
\end{equation}
Substituting \eqref{fay}, the density $\dfb(x)=2n(1-2x)^{n-1}$, which we derive by differentiating \eqref{eq:densB}, and \eqref{eq:condA} into \eqref{eq:condDens}, we have for $0 \leqslant y \leqslant 1$
\begin{equation*}
f_{B|A=y}(x)=\begin{cases}
\frac{2n(n-1)}{f_A(y)}\left(y^{n-2}+((y-2x)^+)^{n-2}-2((2y-1)^+)^{n-2}\right), &0\leqslant x \leqslant \min\{\frac{1-y}{2},\frac{1}{4}\},\\
\frac{1}{f_A(y)}(2n(n-1)y^{n-2}), &\frac{1}{4}\leqslant x \leqslant \frac{1-y}{2}.
\end{cases}
\end{equation*}
The last range is obtained because we formally need that $\frac{1}{4}\leqslant x\leqslant \frac{1}{2}$ and that $0\leqslant y \leqslant 1-2x$, or otherwise that $\frac{1}{4}\leqslant x \leqslant \min\{\frac{1}{2}, \frac{1-y}{2}\}$. However, for $x\in[\frac{1}{4},\frac{1}{2}]$, we have that $y\in[0,\frac{1}{2}]$ and thus the minimum is always achieved at $\frac{1-y}{2}$. We naturally understand the density to be zero if the corresponding interval results to the right limit being smaller than the left limit.

Thus, for $0\leqslant x \leqslant \min\{\frac{1-y}{2},\frac{1}{4}\}$ we obtain
\begin{multline*}
\p[B\leqslant x|A=y]=\int_0^x f_{B|A=y}(z)dz=\\
=\frac{1}{\dfa(y)}\left(2n(n-1)y^{n-2}x-n((y-2x)^+)^{n-1}+ny^{n-1}-4n(n-1)((2y-1)^+)^{n-2}x\right),
\end{multline*}
and for $\frac{1}{4}\leqslant x\leqslant\frac{1-y}{2}$:
\begin{align*}
\p[B\leqslant x|A=y]&=\p[B\leqslant \frac{1}{4}|A=y]+\int_{\frac{1}{4}}^x\frac{1}{f_A(y)}(2n(n-1)y^{n-2})dz\\
&=\frac{1}{f_A(y)}\left(2n(n-1)y^{n-2}x+ny^{n-1}\right).
\end{align*}

\section{Proof of Lemma~\ref{lem:soj}}\label{app:4}
From \eqref{eq2}, we have for $0\leqslant x \leqslant \frac{1}{2}$
\begin{align*}
\Dfs(x)&=\int_0^x\int_0^1 \p[B\leqslant x-y+z|A=y]\dfa(y)\dfs(z)dzdy\\
&=\int_0^x\int_0^{\frac{1+y}{2}-x}\p[B\leqslant x-y+z|A=y]\dfa(y)\dfs(z)dzdy+\int_0^x\int_{\frac{1+y}{2}-x}^1\dfa(y)\dfs(z)dzdy\\
&=\int_0^{(x-\frac{1}{4})^+}\int_0^{\frac{1+y}{2}-x}\p[B\leqslant x-y+z|A=y]\dfa(y)\dfs(z)dzdy\\
&\quad+\int_{(x-\frac{1}{4})^+}^x\int_0^{\frac{1}{4}+y-x}\p[B\leqslant x-y+z|A=y]\dfa(y)\dfs(z)dzdy\\
&\quad+\int_{(x-\frac{1}{4})^+}^x\int_{\frac{1}{4}+y-x}^{\frac{1+y}{2}-x}\p[B\leqslant x-y+z|A=y]\dfa(y)\dfs(z)dzdy\\
&\quad+\int_0^x\int_{\frac{1+y}{2}-x}^1\dfa(y)\dfs(z)dzdy.
\end{align*}
Observe that for the first and the third integral, we have that $\frac{1}{4}\leqslant x-y+z\leqslant\frac{1-y}{2}$. Also observe that for the second integral we have that $0\leqslant x-y+z\leqslant \frac{1}{4}$, and recall \eqref{eq:condB2}. Then we see that the term containing $(2y-1)^+$ is eliminated as $y\leqslant x$ and $0\leqslant x\leqslant\frac{1}{2}$. Summarising these results we have that
\begin{align*}
\Dfs(x)&=\int_0^{(x-\frac{1}{4})^+}\int_0^{\frac{1+y}{2}-x}(2n(n-1)y^{n-2}(x-y+z)+ny^{n-1})\dfs(z)dzdy\\
&\quad+\int_{(x-\frac{1}{4})^+}^x\int_0^{\frac{1}{4}+y-x}(2n(n-1)y^{n-2}(x-y+z)+ny^{n-1})\dfs(z)dzdy\\
&\quad+\int_{(x-\frac{1}{4})^+}^x\int_0^{\frac{1}{4}+y-x}(-n[(y-2(x-y+z))^+]^{n-1})\dfs(z)dzdy\\
&\quad+\int_{(x-\frac{1}{4})^+}^x\int_{\frac{1}{4}+y-x}^{\frac{1+y}{2}-x}(2n(n-1)y^{n-2}(x-y+z)+ny^{n-1})\dfs(z)dzdy\\
&\quad+\int_0^x\int_{\frac{1+y}{2}-x}^1\dfa(y)\dfs(z)dzdy.
\end{align*}
Aggregate the first, second and fourth lines and integrate with respect to $z$ to obtain
\begin{align*}
\Dfs(x)&=\int_0^x(2n(n-1)y^{n-2}(x-y)+ny^{n-1})\Dfs(\tfrac{1+y}{2}-x)dy\\
&\quad+\int_0^x(\frac{1+y}{2}-x)(2n(n-1)y^{n-2})\Dfs(\tfrac{1+y}{2}-x)dy\\
&\quad-\int_0^x\int_0^{\frac{1+y}{2}-x}(2n(n-1)y^{n-2})\Dfs(z)dzdy\\
&\quad-\int_{(x-\frac{1}{4})^+}^x\int_0^{\frac{1}{4}+y-x}n((3y-2x-2z)^+)^{n-1}\dfs(z)dzdy+\int_0^x\dfa(y)(1-\Dfs(\tfrac{1+y}{2}-x))dy,
\end{align*}
which, with further computations, leads to
\begin{align*}
\Dfs(x)&=\int_0^x(n(n-1)y^{n-2}(1-y)+ny^{n-1})\Dfs(\tfrac{1+y}{2}-x)dy\\
&\quad-\int_0^x\int_0^{\frac{1+y}{2}-x}(2n(n-1)y^{n-2})\Dfs(z)dzdy\\
&\quad-\int_{(x-\frac{1}{4})^+}^x n((y-\frac{1}{2})^+)^{n-1}\Dfs(\tfrac{1}{4}+y-x)dy\\
&\quad-\int_{(x-\frac{1}{4})^+}^x\int_0^{\frac{1}{4}+y-x}2n(n-1)((3y-2x-2z)^+)^{n-2}\Dfs(z)dzdy\\
&\quad+\Dfa(x)-\int_0^x\dfa(y)\Dfs(\tfrac{1+y}{2}-x)dy.
\end{align*}
Since $0\leqslant y\leqslant x\leqslant \frac{1}{2}$, the third line integrates to zero. Moreover, the first integral is equal to the last integral, since for $y\leqslant \frac{1}{2}$, $\dfa(y)$ is equal to $(n(n-1)y^{n-2}(1-y)+ny^{n-1})$, cf.\ \eqref{fay}. Thus, for $0\leqslant x\leqslant \frac{1}{2}$, we have that
\begin{multline}\label{eq4}
\Dfs(x)=\Dfa(x)-\int_0^x\int_0^{\frac{1+y}{2}-x}(2n(n-1)y^{n-2})\Dfs(z)dzdy\\
-\int_{(x-\frac{1}{4})^+}^x \int_0^{\frac{1}{4}+y-x} 2n(n-1)((3y-2x-2z)^+)^{n-2}\Dfs(z)dzdy.
\end{multline}

We now examine the case where $\frac{1}{2}\leqslant x\leqslant \frac{3}{4}$. The computations are long, but thoroughly standard. Here we only give explicitly only the most important steps, and indicate all other intermediate computations. Observe that, for this case, we have that $0\leqslant 2x-1\leqslant x-\frac{1}{4}$ and keep in mind Lemma~\ref{lem:condProb}. From \eqref{eq2} we have that
\begin{multline*}
\Dfs(x)=\int_0^{2x-1}\int_0^1 \dfa(y)\dfs(z)dzdy+\int_{2x-1}^x\int_0^{\frac{1+y}{2}-x}\p[B\leqslant x-y+z|A=y]\dfa(y)\dfs(z)dzdy\\+\int_{2x-1}^x\int_{\frac{1+y}{2}-x}^1\dfa(y)\dfs(z)dzdy
\end{multline*}
The second term at the right hand side can be split in four different terms, based on the expressions in Lemma~\ref{lem:condProb}. We examine separately the cases when $y \in (2x-1, x-\frac14)$, $y \in (\frac12, x)$, and when $y \in (x-\frac14, \frac12)$, we examine separately the integrals for $z\in(0,\frac14+y-x)$ and for $z\in(\frac14+y-x, \frac{1+y}{2}-x)$. The resulting six double integrals can be simplified (based on the intervals we work with, some terms are equal to zero), aggregated, and in cases computed (e.g.\ in the first integral in the expression above, the integration with respect to $z$ yields trivially 1), leading to
\begin{align*}
\Dfs(x)&=\int_{2x-1}^x\left(n(n-1)y^{n-2}(1-y)+ny^{n-1}\right)\Dfs(\tfrac{1+y}{2}-x)dy\\
&\quad-\int_{2x-1}^x\int_0^{\frac{1+y}{2}-x}2n(n-1)y^{n-2}\Dfs(z)dzdy\\
&\quad+\int_{x-\frac{1}{4}}^{\frac{1}{2}}\left(-n((y-\frac{1}{2})^+)^{n-1}\right)\Dfs(\tfrac{1}{4}+y-x)dy\\
&\quad-\int_{x-\frac{1}{4}}^{\frac{1}{2}}\int_0^{\frac{1}{4}+y-x}\left(2n(n-1)((3y-2x-2z)^+)^{n-2}\right)\Dfs(z)dzdy\\
&\quad-\int_{\frac{1}{2}}^x\left(n((2y-1)^+)^{n-1}+2n(n-1)(2y-1)^{n-2}(1-y)\right)\Dfs(\tfrac{1+y}{2}-x)dy\\
&\quad+\int_{\frac{1}{2}}^x\int_0^{\frac{1+y}{2}-x}\left(-2n(n-1)((3y-2x-2z)^+)^{n-2}+4n(n-1)(2y-1)^{n-2}\right)\Dfs(z)dzdy\\
&\quad+\Dfa(x)-\int_{2x-1}^x\dfa(y)\Dfs(\tfrac{1+y}{2}-x)dy.
\end{align*}
We observe that the density \dfa, as given by \eqref{fay}, results to having the first and fifth line cancel the last term. Thus, all terms including the quantity $\Dfs(\tfrac{1+y}{2}-x)$ are removed. Moreover, the third line is also trivially equal to zero, which leads to the following integral equation for \Dfs\ when $x\in[\frac12,\frac34]$:
\begin{align}\label{eq9}
\Dfs(x)&=\Dfa(x)-\int_{2x-1}^x\int_0^{\frac{1+y}{2}-x}2n(n-1)y^{n-2}\Dfs(z)dzdy\nonumber\\
&\quad-\int_{x-\frac{1}{4}}^{\frac{1}{2}}\int_0^{\frac{1}{4}+y-x}\left(2n(n-1)((3y-2x-2z)^+)^{n-2}\right)\Dfs(z)dzdy\nonumber\\
&\quad-\int_{\frac{1}{2}}^x\int_0^{\frac{1+y}{2}-x}\left(2n(n-1)((3y-2x-2z)^+)^{n-2}-4n(n-1)(2y-1)^{n-2}\right)\Dfs(z)dzdy.
\end{align}
Observe that, as expected, $\Dfs(x)$ is continuous at $x=\frac12$.

For the last interval, i.e.\ for $\frac{3}{4}\leqslant x\leqslant 1$, we have that $2x-1\geqslant \frac{1}{2}$ and that when $y\geqslant \frac{1}{2}$ we have that $\frac{1}{4}+y-x\geqslant \frac{1+y}{2}-x$. Starting from \eqref{eq2} similar computations that necessitate the usage of Lemma~\ref{lem:condProb} and the density \dfa, given in \eqref{fay}, result to
\begin{align}\label{eq:LastInterval}
\Dfs(x)&=\Dfa(x)-\int_{2x-1}^x\int_0^{\frac{1+y}{2}-x}\left(2n(n-1)y^{n-2}+2n(n-1)((3y-2x-2z)^+)^{n-2}\right)\Dfs(z)dzdy\nonumber\\
&\quad+\int_{2x-1}^x\int_0^{\frac{1+y}{2}-x}\left(4n(n-1)((2y-1)^+)^{n-2}\right)\Dfs(z)dzdy.
\end{align}
As a sanity check, we observe that \Dfs\ is continuous at $x=\frac34$ and that $\Dfs(1)=\Dfa(1)=1$. Aggregating now \eqref{eq4}--\eqref{eq:LastInterval}, we obtain the integral equation of Lemma~\ref{lem:soj}.

\section{Proof of Theorem~\ref{th:3}}\label{app:5}
Let $F_1(x)$ and $F_2(x)$ be two continuous functions in $\mathcal{L}$ with $|F_1-F_2|\leqslant \delta$ for all $x\in [0,1]$. We then have
\begin{align*}
|\Omega F_1(x)-\Omega F_2(x)|&\leqslant\delta\int_{(2x-1)^+}^x\int_0^{\frac{1+y}{2}-x}2n(n-1)y^{n-2}dzdy\\
&+\delta\int_{\max\{(x-\frac{1}{4})^+,2x-1\}}^{x}\int_0^{\frac{1+y}{2}-x}2n(n-1)((3y-2x-2z)^+)^{n-2}dzdy\\
&+\delta\int_{\max\{\frac{1}{2},2x-1\}}^{\max\{\frac{1}{2},x\}}\int_0^{\frac{1+y}{2}-x}4n(n-1)(2y-1)^{n-2}dzdy\\
&=\delta\int_{(2x-1)^+}^x n(n-1)y^{n-2}(1+y)-2n(n-1)y^{n-2}xdy\\
&+\delta\int_{\max\{(x-\frac{1}{4})^+,2x-1\}}^{x} n((3y-2x)^+)^{n-1}-n((2y-1)^+)^{n-1} dy\\
&+\delta\int_{\max\{\frac{1}{2},2x-1\}}^{\max\{\frac{1}{2},x\}}2n(n-1)(2y-1)^{n-2}(1+y-2x)dy.
\end{align*}
Bounding $1+y-2x$ in the last integral by $1-x$ and performing the integration, we have that
\begin{align*}
|\Omega F_1(x)-\Omega F_2(x)|&\leqslant\delta\left[y^{n-1}(n+(n-1)y)-2ny^{n-1}x\right]_{(2x-1)^+}^x\\
&+\delta\left[\frac{1}{3}((3y-2x)^+)^n-\frac{1}{2}((2y-1)^+)^n\right]_{\max\{(x-\frac{1}{4},)^+,2x-1\}}^x\\
&+\delta(1-x)\left[n(2y-1)^{n-1}\right]_{\max\{\frac{1}{2},2x-1\}}^{\max\{\frac{1}{2},x\}}.
\end{align*}
It can be easily shown that the second line can be bounded by ${\delta x^n}/{3}$, thus yielding
\begin{multline*}
|\Omega F_1(x)-\Omega F_2(x)|\leqslant \delta\Big(x^{n-1}\big(n+(n-1)x\big)-2nx^{n}-((2x-1)^+)^{n-1}\big(n+(n-1)(2x-1)^+\big)\Big)\\
+\delta \Big(2nx((2x-1)^+)^{n-1}+\frac{x^n}{3} \Big)+\delta(1-x)\Big(n((2x-1)^+)^{n-1}-n((4x-3)^+)^{n-1}\Big)\\
=\delta\Big(x^{n-1}(n-x-nx)+((2x-1)^+)^{n}\Big)+\delta\Big(\frac{x^n}{3}+(1-x)\big(n((2x-1)^+)^{n-1}-n((4x-3)^+)^{n-1}\big)\Big).
\end{multline*}
Now set
\begin{align*}
f_n(x)&=x^{n-1}(n-x-nx)+((2x-1)^+)^{n}
\intertext{and}
g_n(x)&=(1-x)\left(n((2x-1)^+)^{n-1}-n((4x-3)^+)^{n-1}\right).
\end{align*}
It readily follows that for $x\in[0,1]$ and $n=2$
$$
f_2(x)=2x(1-x)-x^2+((2x-1)^+)^2\leqslant\frac{1}{3}
$$
and that
$$
g_2(x)=(1-x)\left(2(2x-1)^+-2(4x-3)^+\leqslant\frac{1}{4}\right).
$$
Moreover, observe that $\sup_{x\in[0,1]}f_{n+1}(x)\leqslant \sup_{x\in [0,1]}f_{n}(x)$ and likewise for $g_n$. Trivially, $\frac{x^n}{3}\leqslant \frac13$ for all $n>1$ and $x\in[0,1]$.
Thus, we have that
\begin{equation*}
|\Omega F_1(x)-\Omega F_2(x)|\leqslant
\delta\left(\frac{1}{4}+\frac{1}{3}+\frac{1}{3}\right)=\frac{11}{12}\delta<\delta,
\end{equation*}
which implies that $\Omega$ is indeed a contraction mapping on $\mathcal{L}$ with contraction constant $11/12$.

\end{document}